\documentstyle[11pt]{article}

\newcommand{\A}{ {\cal A} }
\newcommand{\B}{ {\cal B} }
\newcommand{\C}{ {\bf C} }
\newcommand{\D}{ {\cal D} }

\newcommand{\M}{ {\cal M} }
\newcommand{\N}{ {\cal N} }

\newcommand{\X}{ {\cal X} }

\newcommand{\Z}{ {\bf Z} }

\newcommand{\ncps}{ ( {\cal A} , \varphi ) }
\newcommand{\Kr}{ \mbox{Kr} } 
\newcommand{\coef}{ \mbox{coef} } 
\newcommand{\perm}{ \mbox{perm} } 
\newcommand{\ins}{ \mbox{ins} }

\newcommand{\freestar}{ \framebox[7pt]{$\star$} }
\newcommand{\freestart}{ \widetilde{\freestar} }
\newcommand{\mtens}{ \odot }

\begin{document}

\title{\bf R-cyclic families of matrices in free probability}
\author{Alexandru Nica \thanks{Research supported by a grant of 
NSERC, Canada.}  \\
Department of Pure Mathematics  \\
University of Waterloo \\
Waterloo, Ontario, N2L 3G1, Canada \\ 
(e-mail: anica@math.uwaterloo.ca)
\and Dimitri Shlyakhtenko \thanks{Partially supported by
an NSF Postdoctoral Fellowship.}  \\
Department of Mathematics \\
U.C.L.A. \\
Los Angeles, CA 90095-1555, USA \\
(email: shlyakht@math.ucla.edu)
\and Roland Speicher \\
Department of Mathematics and Statistics  \\
Queen's University   \\
Kingston, Ontario K7L 3N6, Canada \\
(email: speicher@mast.queensu.ca) }
\date{ }

\maketitle

\vspace{1in}

\begin{abstract}
We introduce the
concept of ``R-cyclic family'' of matrices with entries in a 
non-commutative probability space; the definition consists in 
asking that only the ``cyclic'' non-crossing cumulants of the 
entries of the matrices are allowed to be non-zero.

Let $A_{1}, \ldots , A_{s}$ be an R-cyclic family of $d \times d$
matrices over a non-commutative probability space $\ncps$. We prove 
a convolution-type formula for the explicit computation of the joint
distribution of $A_{1}, \ldots , A_{s}$ (considered in $M_{d} ( \A )$
with the natural state), in terms of the joint distribution 
(considered in the original space $\ncps$) of the entries of the
$s$ matrices.
Several important situations of families of matrices with tractable 
joint distributions arise by application of this formula.

Moreover, let $A_{1}, \ldots , A_{s}$ be a family of $d \times d$
matrices over a non-commutative probability space $\ncps$, let 
$\D \subset M_{d} ( \A )$ denote the algebra of scalar diagonal 
matrices, and let ${\cal C}$ be the subalgebra of $M_{d} ( \A )$
generated by $\{ A_{1}, \ldots , A_{s} \} \cup \D$. We prove that
the R-cyclicity of $A_{1}, \ldots , A_{s}$ is equivalent to a 
property of ${\cal C}$ -- namely that ${\cal C}$ is free from 
$M_{d} ( \C )$, with amalgamation over $\D$.
\end{abstract}

\newpage

\setlength{\baselineskip}{18pt}

{\large\bf Introduction}

$\ $

In the influential paper \cite{V2}, Voiculescu introduced the concepts 
of circular and semicircular systems, and used them to obtain results
about the fundamental groups of the von Neumann algebras associated to 
free groups. There are three main properties of the circulars and 
semicirculars which are essential for the arguments in \cite{V2}:

(a) the compression of a semicircular system by a projection free from 
the system is again a semicircular system;

(b) in the polar decomposition of a circular element, the polar part 
is free from the positive part; 

(c) one can obtain semicircular systems consisting of {\em matrices} 
over a non-commutative probability space, if the entries of these matrices
are themselves chosen (in an appropriate way) to be circular/semicircular
and free.

Each of (a), (b), (c) points to a direction of investigation in the 
combinatorics of free probability.

Concerning (a) and (b), the things are now pretty well understood. 
For (a), we know a general formula describing the 
distribution of the compression by a free projection (see 
\cite{NS1}), or even more generally for what happens when we perform a 
compression by a free matrix unit (see \cite{Sh1}, \cite{N}).
For (b), the relevant class of elements to be studied is the one of 
``R-diagonal elements'', introduced in \cite{NS2}, and which turns out to 
have a lot of good properties (see e.g. \cite{HL}, or
\cite{NSS1}--\cite{NSS3}).

With (c) the situation is not that clear. If we look at the case of only 
one matrix, then the problem is to give effective methods for computing the 
distribution of the matrix, by starting from the joint distribution of 
its entries. Of course, the distribution of the matrix is always 
completely determined by the joint distribution of its entries; the 
issue is here about the word ``effective''. It is unlikely that one can
 give a nice formula which would work in full generality. The problem is 
more like this: to what kind of matrices can one generalize the nice 
facts known about matrices of free circular/semicircular elements? 
We look for a situation which is general enough to contain interesting 
examples, but also particular enough so that a nice formula does exist.

In this paper we propose the concept of R-cyclic matrix (or more
generally, of R-cyclic family of matrices), which we believe is a good
framework for studying the direction (c).

The definition is in terms of the joint R-transform of the entries 
of the matrix -- where the R-transform is the free probabilistic 
counterpart for the characteristic function of the joint distribution.
The coefficients of the R-transform are called non-crossing cumulants. 
The definition of an R-cyclic matrix goes by asking that only
the cyclic non-crossing cumulants of the entries survive; see 
Definition 2.2 in Section 2 below, and see Sections 2.3-2.6 for examples.

If $A$ is an R-cyclic matrix, then all the information about the 
distribution of $A$ is stored in the family of cyclic cumulants of 
its entries. These cyclic cumulants can be in turn nicely stored in
one formal power series $f$ (in $d$ non-commuting variables, where 
$d \times d$ is the size of $A$); the series $f$ is called ``the 
determining series'' of $A$. Our problem is then to find an effective 
method for computing the distribution of the R-cyclic matrix $A$, in terms 
of its determining series $f$. In the Section 2 of the paper we show 
that this problem can be treated by using a convolution-type formula:
\[
(I) \hspace{2cm} R_{A} (z) \ = \ \frac{1}{d} (f \ \freestar \ H_{d} ) 
( \ \underbrace{z, \ldots , z}_{d} \ ),
\]
where: $R_{A}$ is the R-transform of $A$; $H_{d}$ is a certain universal 
series in $d$ indeterminates; and $\freestar$ is a convolution-type 
operation introduced in \cite{NS1}, which appears to play an important role
in combinatorial free probability (see review in Section 1 below).
The formula $(I)$ can be extended to the case of an R-cyclic family of 
matrices (see Definition 2.9 and Theorem 2.10 in Section 2), and can be 
used to obtain various situations when one gets a family of matrices 
with computable joint distribution. Some applications are presented in 
the Section 3 of the paper.

Section 4 is about operations with matrices in an R-cyclic 
family. It is trivial from the definition that if $A_{1}, \ldots , A_{s}$
is an R-cyclic family (of $d \times d$ matrices over a non-commutative 
probability space $\ncps$), then one can add to $A_{1}, \ldots , A_{s}$:

(a) a linear combination of $A_{1}, \ldots , A_{s}$, or

(b) any scalar diagonal matrix,
\newline
and the enlarged family is still R-cyclic. In Lemma 4.2 we show that
a similar statement is true when one adds to $A_{1}, \ldots , A_{s}$ a 
product of some of the matrices in the family; this comes as a fairly easy
application of a formula for non-crossing cumulants with products for 
entries, which was found in \cite{KS}.

The considerations of Section 4 show that the property of a family of 
matrices
$A_{1}, \ldots , A_{s} \in M_{d} ( \A )$ of being R-cyclic is really a 
property of the algebra ${\cal C}$ generated together by 
$A_{1}, \ldots , A_{s}$ and the set of scalar diagonal matrices. The rest
of the paper is devoted to identifying what this property of ${\cal C}$
exactly is. The result turns out to be the following (Theorem 8.2):
\[
(II) \hspace{2.5cm} 
\left\{  \begin{array}{c}
\mbox{the family}  \\
A_{1}, \ldots , A_{s} \\
\mbox{is R-cyclic}
\end{array}  \right\} \ \Leftrightarrow \
\left\{  \begin{array}{c}
\mbox{${\cal C}$ is free from $M_{d} ( \C )$,}  \\
\mbox{with amalgamation over}                    \\
\mbox{scalar diagonal matrices} 
\end{array}  \right\} ,
\]
where $A_{1}, \ldots , A_{s}$ and ${\cal C}$ are as above, and where the
algebra $M_{d} ( \C )$ of scalar $d \times d$ matrices is identified as a
subalgebra of $M_{d} ( \A )$ in the natural way.

In the paper \cite{NSS3} we had shown that an element $a \in \A$ is 
R-diagonal if and only if the matrix 
$\left[ \begin{array}{cc} 0 & a \\ a^{*} & 0 \end{array} \right] \in
M_{2} ( \A )$ is free from $M_{2} ( \C )$, with amalgamation over scalar 
diagonal matrices. But it is easy to see, directly from the definitions, 
that $a$ is R-diagonal if and only if the matrix
$\left[ \begin{array}{cc} 0 & a \\ a^{*} & 0 \end{array} \right]$ is
R-cyclic. Hence the above equivalence $(II)$ can be viewed as an ample
generalization of the named result from \cite{NSS3}.

The equivalence in $(II)$ is obtained by studying non-crossing 
{\em operator-valued} cumulants, in the sense of \cite{S2}; a few basic 
facts about operator-valued cumulants are reviewed in Section 5, and the 
proof of $(II)$ is shown in Section 8. In between 5 and 8 we have two short 
sections where we derive some explicit formulas (used in Section 8) for 
operator-valued cumulants with respect to the algebra $M_{d} ( \C )$ (in 
Section 6), and with respect to the algebra of scalar diagonal matrices 
(in Section 7).

$\ $

$\ $

\setcounter{section}{1}
{\large\bf 1. Basic concepts for the combinatorics of free probability}

$\ $

As a preparation for the theorems proved in Section 2, we review here
a few basic concepts and facts used in combinatorial free probability.
We use the framework of a {\em non-commutative probability 
space}, by which we will simply understand a pair $\ncps$ where $\A$ is 
a complex unital algebra (``the algebra of random variables'') and 
$\varphi : \A \rightarrow \C$ (``the expectation'') is a linear
functional, normalized by $\varphi (1) =1.$ We assume that the reader 
has some familiarity with the concept of freeness for families of 
elements in $\ncps$ (see e.g. \cite{VDN}, Chapter 2).

In the combinatorial study of freeness, an important role is played
by the concepts of {\em moment series} and {\em R-transform} of a 
family of non-commuting random variables. The definition of the first
of these two concepts is straightforward: if $\ncps$ is a non-commutative
probability space, and if $a_{1}, \ldots , a_{s}$ are in $\A$, then 
the numbers of the form:
\begin{equation}
\varphi ( \ a_{r_{1}} \cdots \ a_{r_{n}} \ ), \ \ 
n \geq 1, \ 1 \leq r_{1} , \ldots , r_{n} \leq s,
\end{equation}
are called the {\em joint moments} of $a_{1}, \ldots , a_{s}$; the 
moment series of $a_{1}, \ldots , a_{s}$ 
%denoted $M_{a_{1}, \ldots , a_{s}}$,
is the power series in $s$ non-commuting 
indeterminates $z_{1}, \ldots , z_{s}$ which has the joint moments as
coefficients. That is:
\begin{equation}
M_{a_{1}, \ldots , a_{s}} (z_{1}, \ldots , z_{s}) \ := \ 
\sum_{n=1}^{\infty} \ \sum_{r_{1}, \ldots , r_{n}=1}^{s} \
\varphi ( a_{r_{1}} \cdots a_{r_{n}} ) z_{r_{1}} \cdots z_{r_{n}}.
\end{equation}

The (less straightforward) definition of the R-transform can be placed 
within the framework of a certain convolution operation on formal power 
series which will be used in Section 2, and is reviewed next (in Sections
1.1-1.2, followed by the definition of the R-transform in Section 1.3).

$\ $

{\bf 1.1 Non-crossing partitions.}
Let $\pi = \{ B_{1} , \ldots , B_{k} \}$ be a partition of 
$\{ 1, \ldots , n \}$ -- i.e. $B_{1} , \ldots , B_{k}$ are pairwisely
disjoint non-void sets (called the {\em blocks} of $\pi$), and
$B_{1} \cup \cdots \cup B_{k}$ = $\{ 1, \ldots , n \}$.
We say that $\pi$ is {\em non-crossing} 
if for every $1 \leq i < j < k < l \leq n$ such that $i$ is in the same
block with $k$ and $j$ is in the same block with $l$, it necessarily
follows that all of $i,j,k,l$ are in the same block of $\pi$.
The set of non-crossing partitions of $\{ 1, \ldots , n \}$ will be
denoted by $NC(n).$

For $\pi , \rho \in NC(n),$ we write ``$\pi \leq \rho$'' if each block of
$\rho$ is a union of blocks of $\pi$. Then ``$\leq$'' is a partial order
relation on $NC(n)$, called the refinement order. It turns out that 
$(NC(n), \leq )$ is in fact a lattice, i.e. every two partitions in $NC(n)$
have a lowest upper bound and a greatest lower bound with respect to 
$\leq$.

For $\pi \in NC(n)$ we will denote by $\perm_{\pi}$ the permutation of 
$\{ 1, \ldots , n \}$ which has the blocks of $\pi$ as cycles, in such 
a way that if $B = \{ k_{1} < \cdots < k_{p-1} < k_{p} \}$ is a block 
of $\pi$ then we have
\[
\perm_{\pi} ( k_{1} ) = k_{2}, \ldots ,
\perm_{\pi} ( k_{p-1} ) = k_{p}, 
\perm_{\pi} ( k_{p} ) = k_{1}.
\]
( For example, if $\pi = \{ \ \{ 1,2,5 \} , \{ 3,4 \} \ \} \in NC(5)$, then 
$\perm_{\pi} = \left( \begin{array}{ccccc}
1 & 2 & 3 & 4 & 5  \\
2 & 5 & 4 & 3 & 1
\end{array}  \right)$. )
The set $\{ \perm_{\pi} \ | \ \pi \in NC(n) \}$ has a nice interpretation
in terms of the geometry of the Cayley graph of the symmetric group (see
\cite{B}), and can be a useful instrument in considerations about the 
lattice $NC(n)$.

Unlike the lattice of all partitions of $\{ 1,2, \ldots , n \}$,
$NC(n)$ is anti-isomorphic to itself. We will in fact make extensive 
use of a canonical anti-isomorphism $\Kr: NC(n) \rightarrow NC(n)$,
introduced in \cite{K} and called the Kreweras complementation map.
The map $\Kr$ can be conveniently described by using the permutations
associated to non-crossing partitions, via the following formula:
\begin{equation}
\perm_{\pi} \circ \perm_{\Kr ( \pi )} \ = \ \gamma_{n}, \ \ 
\forall \ \pi \in NC(n),
\end{equation}
where $\gamma_{n}$ is the forward cycle on $\{ 1, \ldots , n \}$
$( \gamma_{n} (1) = 2, \ldots , \gamma_{n} (n-1) =n,
\gamma_{n} (n) = 1).$

$\ $

{\bf 1.2 The operation of boxed convolution.}
Let $s$ be a positive integer. We denote by $\Theta_{s}$ the set of 
all series of the form appearing in Equation (1.2):
\begin{equation}
\Theta_{s} \ = \ \left\{ f \ \begin{array}{ll}
| & f(z_{1}, \ldots , z_{s}) = \sum_{n=1}^{\infty} \
\sum_{r_{1}, \ldots r_{n} = 1}^{s} \ 
\alpha_{r_{1}, \ldots , r_{n}} z_{r_{1}} \cdots z_{r_{n}}   \\
| & \mbox{where } \alpha_{r_{1}, \ldots , r_{n}} \in \C \
(n \geq 1, \ 1 \leq r_{1}, \ldots , r_{n} \leq s)
\end{array}   \right\}  .
\end{equation}
For a series $f$ as in (1.4), we will use the notation 
\begin{equation}
[ \ \mbox{coef } ( r_{1}, \ldots , r_{n} ) \ ] (f)
\end{equation}
to denote the coefficient $\alpha_{r_{1}, \ldots , r_{n}}$ of 
$z_{r_{1}} \cdots z_{r_{n}}$ in $f$.

The operation of boxed convolution, $\freestar$, is an associative
binary operation on the set $\Theta_{s}$. Its definition is inspired 
from the combinatorial theory of convolution in a lattice, as 
developed by Rota and his collaborators (see e.g. \cite{DRS}; the 
lattices of relevance for the definition of $\freestar$ are those 
of non-crossing partitions, $NC(n)$ for $n \geq 1$). 

In order to state the definition of $\freestar$, it is convenient
to first expand the notations for coefficients introduced in (1.5).
If $n \geq 1$, $1 \leq r_{1}, \ldots , r_{n} \leq s$, and if
$B = \{ k_1 < k_2 <  \cdots < k_p \}$ is a non-void subset of 
$\{ 1, \ldots , n \}$, then by ``$( r_{1}, \ldots , r_{n} ) | B$''
we will understand the $p$-tuple 
$(r_{k_{1}} , r_{k_{2}} ,  \ldots , r_{k_{p}} )$ (for example
$(r_1 , r_2 , r_3 , r_4 , r_5 ) | \{ 2,3,5 \} = (r_2 , r_3 , r_5 ))$.
Then for a series $f \in \Theta_{s}$ we introduce the following 
``generalized coefficients'':
\begin{equation}
[ \mbox{coef } (r_{1}, \ldots , r_{n}); \pi ] (f) \ := \ 
\prod_{B \ block \ of \ \pi} \ [ \ \mbox{coef } 
(r_{1}, \ldots , r_{n}) |B \ ] (f) ,
\end{equation}
for every $n \geq 1$, $1 \leq r_{1}, \ldots , r_{n} \leq s$, and for 
every $\pi \in NC(n)$. (For example if $n=4$ and 
$\pi = \{ \{ 1,3 \} , \{ 2 \} , \{ 4 \} \}$, then 
\[
[ \mbox{coef } (r_1 , r_2 , r_3 , r_4 ) ; \pi ] (f) \ = \
[ \mbox{coef } (r_1 , r_3 ) ] (f) \cdot 
[ \mbox{coef } (r_2 ) ] (f) \cdot
[ \mbox{coef } (r_4 ) ] (f) ,
\]
for any $1 \leq r_1 , r_2 , r_3 , r_4 \leq s.)$

By using the notation introduced in (1.6), the boxed convolution 
$f \ \freestar \ g$ of two series $f,g \in \Theta_{s}$
is described by the formula:
\begin{equation}
[ \mbox{coef }(r_{1}, \ldots ,r_{n}) ] (f \ \freestar \ g) \ := \
\end{equation}
\[
\sum_{\pi \in NC(n)}
[ \mbox{coef }(r_{1}, \ldots , r_{n}); \pi ] (f) \cdot 
[ \mbox{coef }(r_{1}, \ldots , r_{n}); \Kr( \pi ) ] (g) ,
\]
holding for every $n \geq 1$ and $1 \leq r_{1}, \ldots , r_{n} \leq s$,
and where $\Kr ( \pi )$ is the Kreweras complement of the 
partition $\pi \in NC(n)$.

It can be shown that $\freestar$ is associative and unital, where the 
unit is the series $\Delta (z_{1}, \ldots , z_{s} )$
$:= z_{1} + \cdots + z_{s}$.  A series $f \in \Theta_{s}$ is invertible 
with respect to $\freestar$ if and only if its coefficients of degree 1,
$[ \mbox{coef } (r) ] (f)$, $1 \leq r \leq s$, are all different
from 0 (see \cite{NS1}, Section 3).

$\ $

{\bf 1.3 R-transform and free cumulants.} Let $a_{1}, \ldots , a_{s}$ be 
an $s$-tuple of elements in a non-commutative probability space $\ncps$. 
The R-transform of the $s$-tuple, $R_{a_{1}, \ldots , a_{s}}$, is a 
series in the set $\Theta_{s}$ of Equation (1.4). A succinct
way of introducing $R_{a_{1}, \ldots , a_{s}}$ goes by using the 
boxed convolution $\freestar$ and a special series 
$\mbox{M\"ob}_{s} \in \Theta_{s}$, called the M\"obius series.

$\mbox{M\"ob}_{s}$ is defined as the inverse under $\freestar$ of the 
``zeta series in $s$ indeterminates'',
\begin{equation}
\mbox{Zeta}_{s} ( z_{1}, \ldots , z_{s} ) \ := \ 
\sum_{n=1}^{\infty} \ \sum_{ r_{1}, \ldots , r_{n} =1}^{s} \
z_{r_{1}} \cdots z_{r_{n}}.
\end{equation}
It is not hard to determine the coefficients of $\mbox{M\"ob}_{s}$ 
explicitly:
\begin{equation}
\mbox{M\"ob}_{s} ( z_{1}, \ldots , z_{s} ) \ = \ 
\sum_{n=1}^{\infty} \ \sum_{r_{1}, \ldots , r_{n}=1}^{s}  
(-1)^{n+1}  \frac{ (2n-2)! }{ (n-1)! n! } \ 
z_{r_{1}} \cdots  z_{r_{n}}
\end{equation}
(see e.g. \cite{NS1}, Remark 3.8).

Now, if $\ncps$ is a non-commutative probability space, and if 
$a_{1}, \ldots , a_{s} \in \A$, then we define:
\begin{equation}
R_{a_{1}, \ldots , a_{s}} \ := \ M_{a_{1}, \ldots , a_{s}} \
\freestar \ \mbox{M\"ob}_{s},
\end{equation}
where $M_{a_{1}, \ldots , a_{s}}$ is the moment series from 
Equation (1.2). It is clear that $R_{a_{1}, \ldots , a_{s}}$ 
contains the same information about $a_{1}, \ldots , a_{s}$ as 
the moment series, since Equation (1.10) can be re-written 
equivalently as
\begin{equation}
M_{a_{1}, \ldots , a_{s}} \ = \ R_{a_{1}, \ldots , a_{s}} \
\freestar \ \mbox{Zeta}_{s}.
\end{equation}

Following \cite{S1}, it is customary to denote the coefficient of 
$z_{r_{1}} \cdots z_{r_{n}}$ in $R_{a_{1}, \ldots , a_{s}}$ by:
\begin{equation}
k_{n} ( a_{r_{1}}, \ldots , a_{r_{n}} ).
\end{equation}
More generally, given $n \geq 1$, $1 \leq r_{1}, \ldots , r_{n} \leq s$,
and a partition $\pi \in NC(n)$, we use the notation
\begin{equation}
k_{\pi} ( a_{r_{1}}, \ldots , a_{r_{n}} )
\end{equation}
for the ``generalized coefficient''
$[ \mbox{coef } (r_{1}, \ldots , r_{n}); \pi ](R_{a_{1}, \ldots ,a_{s}})$
defined as in Equation (1.6). These generalized coefficients are called 
the {\em non-crossing cumulants} of the $s$-tuple $a_{1}, \ldots , a_{s}$.
It is worth keeping in mind that for any $n \geq 1$ and $\pi \in NC(n)$, 
it makes sense to view $k_{\pi}$ as a multilinear map from $\A^{n}$ to 
$\C$ (see \cite{S1}).

$\ $

{\bf 1.4 R-transform and freeness.}
The R-transform and the boxed convolution turn out to have very pleasant
properties in connection to the addition and multiplication 
of free $n$-tuples -- see \cite{V1}, \cite{NS1}. Even more 
importantly, R-transforms (or equivalently, non-crossing cumulants)
can be used to provide a neat characterization of freeness. To be 
precise: let $a_{1}', \ldots , a_{m}', a_{1}'', \ldots , a_{n}''$ be
elements of the non-commutative probability space $\ncps$; then the 
freeness of the families $\{ a_{1}', \ldots , a_{m}' \}$ and 
$\{ a_{1}'', \ldots , a_{n}'' \}$ is equivalent to the equation
\begin{equation}
R_{a_{1} ' , \ldots , a_{m} ' , a_{1} '' , \ldots , a_{n} '' } 
( z_{1}' , \ldots , z_{m} ' , z_{1} '' , \ldots , z_{n} '' ) \ = \ 
\end{equation}
\[
= \  R_{a_{1} ' , \ldots , a_{m} ' } 
( z_{1}' , \ldots , z_{m} ' )  +  
R_{a_{1} '' , \ldots , a_{n} '' } 
( z_{1} '' , \ldots , z_{n} '' ).
\]
It is obvious how Equation (1.14) extends by induction to the case of $s$
(instead of just two) families of elements. Note that in the case of $s$ 
families having one element each, we obtain the following: the elements 
$a_{1}, \ldots , a_{s} \in \A$ form a free family if and only if we have
that
\begin{equation}
R_{a_{1}, \ldots , a_{s}} ( z_{1}, \ldots , z_{s} ) \ = \ 
R_{a_{1}} (z_{1}) + \cdots + R_{a_{s}} ( z_{s} ).
\end{equation}

$\ $

{\bf 1.5 Extended boxed convolution.} Let $s$ and $d$ be positive integers.
Consider the set $\Theta_{sd}$ of power series in $sd$ non-commuting
indeterminates $z_{1,1}, \ldots , z_{r,i}, \ldots , z_{s,d}$. The same
formula as in Equation (1.7) above can be used to define a ``convolution
operation'', denoted in what follows by $\freestart$, which gives a right 
action of $\Theta_{d}$ on $\Theta_{sd}$. More precisely,
if $f \in \Theta_{sd}$ and $g \in \Theta_{d}$ then we define
$f \ \freestart \ g  \in \Theta_{sd}$ by the following formula:
\begin{equation}
[ \mbox{coef }( \ (r_{1},i_{1}), \ldots ,(r_{n},i_{n}) \ ) ] 
(f \ \freestart \ g) \ := 
\end{equation}
\[
\sum_{\pi \in NC(n)}
[ \mbox{coef }( \ (r_{1},i_{1}), \ldots , (r_{n},i_{n}) \ ); \pi ] (f) 
\cdot [ \mbox{coef }(i_{1}, \ldots , i_{n}); \Kr( \pi ) ] (g) ,
\]
holding for every $n \geq 1$ and for every 
$1 \leq r_{1}, \ldots , r_{n} \leq s$, $1 \leq i_{1}, \ldots , i_{n} \leq d$.
Some trivial adjustments of the considerations made in Section 4 of 
\cite{NS1} for $\freestar$ show that $\freestart$ is indeed a right 
action of $\Theta_{d}$ on $\Theta_{sd}$, in the sense that the equation
\begin{equation}
( f \ \freestart \ g ) \ \freestart \ h \ = \ 
f \ \freestart \ ( g \ \freestar \ h )
\end{equation}
holds for every $f \in \Theta_{sd}$ and $g,h \in \Theta_{d}$. 

Let us also record the fact that:
\begin{equation}
f \ \freestart \ \mbox{Zeta}_{d} \ = \ 
f \ \freestar \ \mbox{Zeta}_{sd}, \ \ \ \forall \ f \in \Theta_{sd}
\end{equation}
(where on the right-hand side of (1.18), $\freestar$ denotes the boxed 
convolution operation on $\Theta_{sd}$). This relation is obvious if
one takes into account the fact that 
any Zeta series has all the coefficients equal to 1.

 From (1.17) and (1.18) it is immediate that one also has:
\begin{equation}
f \ \freestart \ \mbox{M\"ob}_{d} \ = \ 
f \ \freestar \ \mbox{M\"ob}_{sd}, \ \ \ \forall \ f \in \Theta_{sd}.
\end{equation}
Note that, as a consequence, we can write the relation
\begin{equation}
M_{a_{1,1}, \ldots , a_{r,i}, \ldots , a_{s,d}} \ \freestart \ 
\mbox{M\"ob}_{d} \ = \ 
R_{a_{1,1}, \ldots , a_{r,i}, \ldots , a_{s,d}},
\end{equation}
holding for any family 
$\{ a_{r,i} \ | \ 1 \leq r \leq s, \ 1 \leq i \leq d \}$
of elements in some non-commutative probability space $\ncps$.

$\ $

{\bf 1.6 Dilations and scalar multiples of power series.} Let $s$ be a 
positive integer, let $f$ be a series in $\Theta_{s}$, and let $\alpha$ be 
a complex number. We denote by $f \circ D_{\alpha}$ the series in
$\Theta_{s}$ which is defined by the equation:
``$(f \circ D_{\alpha} ) (z_{1}, \ldots , z_{s})$ =
$f ( \alpha z_{1}, \ldots , \alpha z_{s})$'', or more rigorously by the 
fact that:
\[
[ \coef (r_{1}, \ldots , r_{n}) ] (f \circ D_{\alpha}) \ = \ 
\alpha^{n}  \cdot [ \coef (r_{1}, \ldots , r_{n}) ] (f), \ \ 
\forall \ n \geq 1, \ 1 \leq r_{1}, \ldots , r_{n} \leq s.
\]

The formulas relating $\freestar$ with dilation and with scalar 
multiplication which are proved in \cite{NS1} can be easily extended to
the case of $\freestart$. Concerning dilation we have:
\begin{equation}
(f \circ D_{\alpha} ) \ \freestart \ g \ = \ 
f \ \freestart \ (g \circ D_{\alpha}) \ = \ 
(f \ \freestart \ g) \circ D_{\alpha},
\end{equation}
for every $f \in \Theta_{ds}$, $g \in \Theta_{d}$, $\alpha \in \C$.
Concerning scalar multiplication we have the formula:
\begin{equation}
(\alpha f) \ \freestart \ ( \alpha g) \ = \ 
\alpha ( \ ( f \ \freestart \ g) \circ D_{\alpha} ), \ \ \forall \ 
f \in \Theta_{ds}, \ g \in \Theta_{d}, \ \alpha \in \C .
\end{equation}
It is sometimes convenient to use Equation (1.22) in the form:
\begin{equation}
( \alpha f ) \ \freestart \ g \ = \ 
\alpha ( \ f \ \freestart \ ( \frac{1}{\alpha} g \circ D_{\alpha} ) \ ),
\end{equation}
holding for $f \in \Theta_{ds}$, $g \in \Theta_{d}$, and
$\alpha \in \C \setminus \{ 0 \}$.
 
$\ $

{\bf 1.7 The special series $H_{d}$.} Let $d$ be a positive integer.
In this paper we also encounter the ``geometric series in $d$ separate
indeterminates'',
\begin{equation}
G_{d} ( z_{1}, \ldots , z_{d}) \
= \ \sum_{n=1}^{\infty} \ \sum_{i=1}^{d}  z_{i}^{n} \ \
( \ = \frac{z_{1}}{1-z_{1}} + \cdots + \frac{z_{d}}{1-z_{d}} \ ),
\end{equation}
and a series derived from $G_{d}$ which can be described as follows:
\begin{equation}
H_{d} \ := \ G_{d} \ \freestar \ (d \cdot \mbox{M\"ob}_{d} \circ D_{1/d}) .
\end{equation}
To give an idea of how $H_{d}$ looks like, here is its truncation to
order three:
\[
H_{d} (z_{1}, \ldots , z_{d}) \ = \ \sum_{i=1}^{d} z_{i} \ + \
\sum_{i_{1},i_{2}=1}^{d} ( \delta_{i_{1},i_{2}} - \frac{1}{d} )
z_{i_{1}}z_{i_{2}} 
\]
\[
+ \ \sum_{i_{1},i_{2},i_{3}=1}^{d} \ \Bigl( \delta_{i_{1},i_{2},i_{3}} -
\frac{1}{d} ( \delta_{i_{1},i_{2}} + \delta_{i_{1},i_{3}} +
\delta_{i_{2},i_{3}} ) + \frac{2}{d^{2}} \ \Bigr)
z_{i_{1}}z_{i_{2}}z_{i_{3}} \ + \cdots  
\]

Note that a direct application of Equation (1.23) (in the particular 
case when $\freestart$ is $\freestar$ on $\Theta_{d}$, and 
$\alpha = 1/d$) gives the alternative formula:
\begin{equation}
H_{d} \ = \ d \cdot \Bigl( \ ( \frac{1}{d} G_{d} ) \ \freestar \
\mbox{M\"ob}_{d} \ \Bigr) .
\end{equation}
Furthermore, the latter equation has the following interpretation.
Let $tr_{d}$ denote the normalized trace on the algebra $M_{d} ( \C )$,
and consider the matrices $P_{1}, \ldots , P_{d} \in M_{d} ( \C )$
where $P_{i}$ has its $(i,i)$-entry equal to 1 and all the other 
entries equal to 0. Then, obviously:
\[
M_{P_{1}, \ldots , P_{d}} \ = \ \frac{1}{d} G_{d}
\]
(moment series considered in the non-commutative probability space
$( M_{d} ( \C ), tr_{d} )$ ); hence:
\[
( \frac{1}{d} G_{d} ) \ \freestar \ \mbox{M\"ob}_{d} \ = \
M_{P_{1}, \ldots , P_{d}} \ \freestar \ \mbox{M\"ob}_{d} \ = \
R_{P_{1}, \ldots , P_{d}},
\]
and the formula (1.26) for $H_{d}$ takes the form
\begin{equation}
H_{d} \ = \ d \cdot R_{P_{1}, \ldots , P_{d}}.
\end{equation}

An application of Equation (1.27) is that for every $n \geq 2$, every 
$k \in \{ 1, \ldots , n \}$, and every fixed indices 
$i_{1}, \ldots , i_{k-1}, i_{k+1}, \ldots , i_{n} \in \{ 1, \ldots d \}$,
we have:
\begin{equation}
\sum_{i=1}^{d} [ \coef (i_{1}, \ldots , i_{k-1}, i, i_{k+1}, \ldots ,
i_{n} ] (H_{d}) = 0.
\end{equation}
Indeed, the sum on the left-hand side of (1.28) is equal to:
\[
d \cdot \sum_{i=1}^{d}
k_{n} ( P_{i_{1}}, \ldots, P_{i_{k-1}}, P_{i}, P_{i_{k+1}}, 
\ldots , P_{i_{n}} ) \ 
\mbox{ (by (1.27))}
\]
\[
= \ k_{n} ( P_{i_{1}}, \ldots, P_{i_{k-1}}, I, P_{i_{k+1}}, 
\ldots , P_{i_{n}} )  \
\mbox{ (by the multilinearity of $k_{n}$)},
\]
and the latter quantity equals 0 by (1.14) and the fact that the identity 
matrix $I$ is free from $\{ P_{1}, \ldots , P_{d} \}$ in 
$( M_{d} ( \C ) , tr_{d} )$.

$\ $

$\ $

\setcounter{section}{2}
{\large\bf 2. R-cyclic matrices and their R-transforms}

$\ $

{\bf 2.1 Notation.} Let $\ncps$ be a non-commutative probability space,
\setcounter{equation}{0}
and let $d$ be a positive integer. Consider the algebra $M_{d} ( \A )$ of 
$d \times d$ matrices over $\A$. We denote by $\varphi_{d}$ the linear 
functional on $M_{d} ( \A )$ defined by the formula:
\begin{equation}
\varphi_{d} ( \ [a_{i,j}]_{i,j=1}^{d} \ ) \ = \ 
\frac{1}{d} \sum_{i=1}^{d} \varphi (a_{i,i}).
\end{equation}
Then $( M_{d} ( \A ), \varphi )$ is a non-commutative probability space, 
too.

$\ $

{\bf 2.2 Definition.} Let $\ncps$ and $d$ be as above. 
A matrix $A = [a_{i,j}]_{i,j=1}^{d} \in M_{d} ( \A )$ is said to be 
{\em R-cyclic} if the following condition holds:
\[
k_{n} ( a_{i_{1},j_{1}}, \ldots , a_{i_{n},j_{n}} ) \ = \ 0
\]
for every $n \geq 1$ and every
$1 \leq i_{1},j_{1}, \ldots , i_{n},j_{n} \leq d$ for which it is not 
true that $j_{1} = i_{2}, \ldots , j_{n-1}=i_{n}, j_{n}= i_{1}$.

If the matrix $A$ is R-cyclic, then the series:
\begin{equation}
f (z_{1}, \ldots , z_{d}) \ := \ \sum_{n=1}^{\infty} \
\sum_{i_{1}, \ldots , i_{n} =1}^{d} \ 
k_{n} (a_{i_{n},i_{1}}, a_{i_{1},i_{2}}, \ldots , a_{i_{n-1},i_{n}} )
z_{i_{1}} z_{i_{2}} \cdots z_{i_{n}}
\end{equation}
is called the {\em determining series} of $A$.

$\ $

{\bf 2.3 Example.} Consider a diagonal matrix,
\[
A \ := \ \left[ \begin{array}{ccc}
a_{1}    &          &   0    \\
         & \ddots   &        \\
0        &          & a_{d}
\end{array}  \right] \in M_{d} ( \A ),
\]
where $\ncps$ and $d$ are as above. An application of Equation (1.15) 
shows that $A$ is R-cyclic if and only if the elements 
$a_{1}, \ldots , a_{d}$ form a free family; if this is the case, then the
determining series of $A$ coincides with the joint R-transform 
$R_{a_{1}, \ldots , a_{d}}$.

$\ $

For more elaborate examples we will use the framework of a 
$*$-probability space, which is also the one most frequently encountered
in applications. A $*$-probability space is a non-commutative probability
space $\ncps$ where $\A$ is a $*$-algebra, and $\varphi$ has the property 
that $\varphi ( a^{*} ) = \overline{\varphi (a)}$, $\forall \ a \in \A$.

$\ $

{\bf 2.4 Example.} Let $\ncps$ be a $*$-probability space, and let 
$\{ e_{i,j} \ | \ 1 \leq i,j \leq d \}$ be a family of elements of $\A$
which satisfy the following relations: $e_{i,j}^{*}=e_{j,i}$ for all 
$1 \leq i,j \leq d$, $e_{i,j}e_{k,l} = \delta_{j,k}e_{i,l}$ for all 
$1 \leq i,j,k,l \leq d$, and $\sum_{i=1}^{d} e_{i,i} = I$. We will 
assume in addition that $\varphi (e_{i,j}) = 0$ whenever $i \neq j$, and
that $\varphi (e_{1,1})= \cdots = \varphi ( e_{d,d} ) = 1/d$. We denote
by $( {\cal C}, \psi )$ the compression of $\ncps$ by $e_{1,1}$, i.e:
\[
{\cal C} := e_{1,1} \A e_{1,1}, \ \ 
\psi := d \cdot \varphi | {\cal C}.
\]

Let now $a$ be a selfadjoint element of $\A$, which is free from 
$\{ e_{i,j} \ | \ 1 \leq i,j \leq d \}$. We compress $a$ by the matrix 
unit formed by the $e_{i,j}$'s, and we move the compressions under the 
projection $e_{1,1}$; that is, we consider the family of elements:
\[
c_{i,j} \ := \ e_{1,i}ae_{j,1} \in {\cal C}, \ \ 
1 \leq i,j \leq d.
\]
One can compute explicitly the free cumulants of the family 
$\{ c_{i,j} \ | \ 1 \leq i,j \leq d \}$, and obtain that for every 
$n \geq 1$ and $1 \leq i_{1}, j_{1}, \ldots , i_{n},j_{n} \leq d$:
\[
k_{n}( c_{i_{1},j_{1}}, \ldots , c_{i_{n},j_{n}} ) \ = \ 
\left\{  \begin{array}{ll}
d^{-(n-1)} k_{n}( a, \ldots , a)  &
\mbox{if } j_{1}=i_{2}, \ldots , j_{n-1}=i_{n}, j_{n}=i_{1}  \\
                                     &                       \\
0                                   & \mbox{otherwise}
\end{array}  \right.
\]
(see Theorem 8.14 or Theorem 17.3 in the notes \cite{NS3}). In other 
words, the matrix $C = [ c_{i,j} ]_{i,j=1}^{d} \in M_{d} ( {\cal C} )$
is R-cyclic, with determining series:
\[
f( z_{1}, \ldots , z_{d} ) \ = \ 
\sum_{n=1}^{\infty} \ \sum_{i_{1}, \ldots , i_{n}=1}^{d} \
d^{-(n-1)} k_{n}(a, \ldots , a) z_{i_{1}} \cdots z_{i_{n}}
\]
\[
= \ d \cdot \sum_{n=1}^{\infty} k_{n}(a, \ldots , a) \cdot 
\Bigl( \frac{z_{1}+ \cdots +z_{d}}{d} \Bigr)^{n}
\]
\[
= \ d \cdot R_{a} ( \ \frac{z_{1}+ \cdots + z_{d}}{d} \ ),
\]
where $R_{a}$ is the R-transform of $a$, in the original space $\ncps$.

$\ $

{\bf 2.5 Example.} Let $\ncps$ be a $*$-probability space. Let $a \in \A$
be an {\em R-diagonal element}, by which we mean 
that the joint R-transform of $a$ and $a^{*}$ is of the form
\[
R_{a,a^{*}} (z_{1},z_{2}) \ = \ \sum_{n=1}^{\infty} \alpha_{n}
( \ (z_{1}z_{2})^{n} + (z_{2}z_{1})^{n} \ )
\]
for a sequence of real coefficients $( \alpha_{n} )_{n=1}^{\infty}$ (see 
\cite{NS2}). The series $f(z) := \sum_{n=1}^{\infty} \alpha_{n} z^{n}$ 
is called the determining series of $a$.

Now consider the non-commutative probability space 
$( M_{2} ( \A ) , \varphi_{2} )$ defined as in Section 2.1, and the
selfadjoint matrix:
\[
A \ = \ \left[ \begin{array}{ll}
0     & a  \\
a^{*} & 0
\end{array}  \right] \ \in \ M_{2} ( \A ).
\]
One immediately checks that $A$ is R-cyclic (and in fact that also 
conversely, the R-cyclicity of $A$ implies the R-diagonality of $a$).
Moreover, the determining series of $A$ (as defined in Section 2.2) 
coincides with the determining series of the R-diagonal element $a$.
A number of results known about R-diagonal elements can be incorporated
in the theory of R-cyclic matrices by using this trick.

$\ $

{\bf 2.6 Example.} The situation discussed in the Example 2.5 can be
generalized to the one of a selfadjoint matrix with free R-diagonal entries.
More precisely, let $\ncps$ be a $*$-probability space, let $d$ be a 
positive integer, and suppose that the elements 
$\{ a_{i,j} \ | \ 1 \leq i,j \leq d \}$ of $\A$ have the following 
properties:

(i) $a_{i,j}^{*} = a_{j,i}$, $\forall \ 1 \leq i,j \leq d$;

(ii) $a_{i,j}$ is R-diagonal whenever $i \neq j$;

(iii) the $d(d+1)/2$ families:
$\{ a_{i,i} \}$ for $1 \leq i \leq d$, together with
$\{ a_{i,j}, a_{j,i} \}$ for $1 \leq i<j \leq d$, are free in $\ncps$.
\newline
Then the matrix $A := [ a_{i,j} ]_{i,j=1}^{d} \in M_{d} ( \A )$ is 
R-cyclic. Indeed, the freeness condition (iii) combined with the 
R-diagonality of $a_{i,j}$ for $i \neq j$ implies that
the only free cumulants made with the entries of $A$ which 
could possibly be non-zero are:
\[
\left\{  \begin{array}{ll}
k_{n} ( a_{i,i}, \ldots , a_{i,i} ) & \mbox{with} \
n \geq 1, \ 1 \leq i \leq d, \mbox{and}              \\
k_{n} ( a_{i,j}, a_{j,i}, \ldots , a_{i,j}, a_{j,i} ) & \mbox{with }
n \geq 1 \mbox{ even}, \ 1 \leq i,j \leq d, \ i \neq j; 
\end{array}  \right.
\]
all these cumulants fall within the pattern allowed by the
definition of R-cyclicity.

$\ $

{\bf 2.7 Remarks.} 1) Variations of the Example 2.6 can be fabricated, 
such that nonselfadjoint matrices are obtained. For this purpose, it is 
more natural to use the concept of R-cyclic family of matrices, given
in Definition 2.9 below, and the $d \times d$ matrix which appears 
should be considered together with its adjoint. We mention that a 
particularly intriguing construction of this type -- upper triangular 
$d \times d$ matrix with circular $*$-distribution -- was studied 
recently in \cite{DH}.

2) In the Example 2.6 one can take the $a_{i,j}$'s to be 
circular/semicircular, thus obtaining a matrix $A$ as considered in 
\cite{V2}. Recall that $a \in \A$ is said to be {\em semicircular} of 
radius $r$ if $a=a^{*}$ and if
\[
\varphi ( a^{n} ) \ = \ \frac{2}{\pi r^{2}} \int_{-r}^{r} 
t^{n} \sqrt{r^{2}-t^{2}} \ dt, \ \ \forall \ n \geq 1;
\]
and that $c \in \A$ is said to be {\em circular} of radius $r$ if it is 
of the form $c = (a+ib)/ \sqrt{2}$, where each of $a,b$ is semicircular 
of radius $r$, and $a$ is free from $b$. It can be shown (see e.g. 
\cite{VDN}, Chapter 3) that if $a \in \A$ is semicircular of radius $r$,
then $k_{2} ( a,a ) = r^{2}/4$ and $k_{n} ( a,a, \ldots , a ) = 0$ for 
$n \neq 2$. As an easy consequence (see e.g. \cite{NS2}), a circular element 
$c$ of radius $r$ is R-diagonal with 
$R_{c,c^{*}} ( z_{1},z_{2} ) = ( r^{2}/4) \cdot (z_{1}z_{2}+z_{2}z_{1})$.
Thus an example of R-cyclic matrix $A = [ a_{i,j} ]_{i,j=1}^{d}$ is provided
by the case when every $a_{i,i}$ is semicircular (of some radius $r_{i,i}$),
every $a_{i,j}$ with $i \neq j$ is circular (of some radius $r_{i,j}$), 
and the conditions (i), (iii) of Example 2.6 are satisfied.

$\ $

The following theorem indicates how the distribution of an R-cyclic
matrix (considered in the non-commutative probability space 
$( M_{d} ( \A ), \varphi_{d} )$ ) can be obtained from the determining 
series of the matrix.

$\ $

{\bf 2.8 Theorem.} Suppose that $A$ is an R-cyclic matrix, and let $f$ 
denote the determining series of $A$. Then we have the formulas:
\begin{equation}
M_{A} (z) \ = \ \frac{1}{d} ( f \ \freestar \ G_{d}) 
( \ \underbrace{z, \ldots , z}_{d \ times} \ ),
\end{equation}
and
\begin{equation}
R_{A} (z) \ = \ \frac{1}{d} ( f \ \freestar \ H_{d})
( \ \underbrace{z, \ldots , z}_{d \ times} \ ),
\end{equation}
where the series $G_{d}$ and $H_{d}$ are as defined in Section 1.7.

$\ $

Before starting on the proof of Theorem 2.8, it is convenient to observe
that the discussion about R-cyclicity can be generalized without much effort
to the situation of a family of matrices, as follows.

$\ $

{\bf 2.9 Definition.} Let $\ncps$ be a non-commutative probability space,
and let $d$ be a positive integer. Let
$A_{1} = [ a_{i,j}^{(1)} ]_{i,j=1}^{d}, \ldots ,
A_{s} = [ a_{i,j}^{(s)} ]_{i,j=1}^{d}$ be matrices in $M_{d} ( \A )$.
We say that the family $A_{1}, \ldots , A_{s}$ is {\em R-cyclic} if 
the following condition holds:
\[
k_{n} ( a_{i_{1},j_{1}}^{(r_{1})}, \ldots ,
a_{i_{n},j_{n}}^{(r_{n})} ) \ = \ 0,
\]
for every $n \geq 1$, every $1 \leq r_{1}, \ldots , r_{n} \leq s$, and
every $1 \leq i_{1},j_{1}, \ldots , i_{n},j_{n} \leq d$ for which it is 
not true that $j_{1} = i_{2}, \ldots , j_{n-1}=i_{n}, j_{n}= i_{1}$.

If the family $A_{1}, \ldots , A_{s}$ is R-cyclic, then the power series 
in $ds$ indeterminates:
\begin{equation}
f (z_{1,1}, \ldots , z_{s,d}) \ := \
\end{equation}
\[
\sum_{n=1}^{\infty} \
\sum_{i_{1}, \ldots , i_{n} =1}^{d} \ 
\sum_{r_{1}, \ldots , r_{n} =1}^{s} \ 
k_{n} ( \ a_{i_{n},i_{1}}^{(r_{1})}, a_{i_{1},i_{2}}^{(r_{2})}, \ldots , 
a_{i_{n-1},i_{n}}^{(r_{n})} \ )
\cdot z_{r_{1},i_{1}} z_{r_{2},i_{2}} \cdots z_{r_{n},i_{n}}
\]
is called the {\em determining series} of the family.

$\ $

{\bf 2.10 Theorem.} Suppose that $A_{1}, \ldots , A_{s}$ is an R-cyclic family
of matrices, with determining series $f$. Then we have the formulas:
\begin{equation}
M_{A_{1}, \ldots , A_{s}} (z_{1}, \ldots , z_{s}) \ = \ \frac{1}{d} 
( f \ \freestart \ G_{d}) ( \
\underbrace{z_{1}, \ldots , z_{1},}_{d \ times} \ldots ,
\underbrace{z_{s}, \ldots , z_{s}}_{d \ times}  \ ),
\end{equation}
and
\begin{equation}
R_{A_{1}, \ldots , A_{s}} (z_{1}, \ldots , z_{s}) \ = \ \frac{1}{d} 
( f \ \freestart \ H_{d}) ( \
\underbrace{z_{1}, \ldots , z_{1},}_{d \ times} \ldots ,
\underbrace{z_{s}, \ldots , z_{s}}_{d \ times}  \ ),
\end{equation}
where the operation $\freestart$ is as described in Section 1.5, and
where the series $G_{d}$ and $H_{d}$ are as defined in Section 1.7.

$\ $

In the proof of Theorem 2.10 we will use the following lemma:

$\ $

{\bf 2.11 Lemma.} Consider the framework of Theorem 2.10. Let $n$ be
a positive integer, let $\pi$ be in $NC(n)$, and consider some 
indices $1 \leq r_{1}, \ldots , r_{n} \leq s$,
$1 \leq i_{1}, \ldots , i_{n} \leq d$. Then we have the equality:
\begin{equation}
k_{\pi} ( \ a_{i_{n},i_{1}}^{(r_{1})},
a_{i_{1},i_{2}}^{(r_{2})}, \ldots , a_{i_{n-1},i_{n}}^{(r_{n})} \ ) \ = 
\end{equation}
\[
[ \coef ( \ (r_{1},i_{1}), \ldots , (r_{n},i_{n}) \ ); \pi \ ] (f)
\cdot [ \coef ( \ i_{1}, \ldots , i_{n} \ ); \Kr( \pi ) \ ] (G_{d}).
\]

$\ $

{\bf Proof.} We will work with the permutations associated to $\pi$ and 
to $Kr ( \pi )$ (as discussed in Section 1.1). We will use cyclic 
notations modulo $n$ for indices -- i.e, ``$i_{k+1}$'' will mean 
``$i_{1}$'' if $k=n$ and ``$i_{k-1}$'' will mean ``$i_{n}$'' if $k=1$.
 
Since every coefficient $G_{d}$ is equal either 
to 0 or to 1, the generalized coefficient of $G_{d}$ appearing 
on the right-hand side of (2.8) also is 0 or 1. So we have two cases.

\vspace{6pt}

{\em Case 1:}
$[ \coef ( \ i_{1}, \ldots , i_{n} \ ); \Kr( \pi ) \ ] (G_{d}) =1.$

\vspace{6pt}

By writing explicitly what the generalized coefficient of $G_{d}$ is, 
we find that:
\begin{equation}
\left\{
\mbox{  \begin{tabular}{c}
$1 \leq k,l \leq n$,                           \\
$k,l$ in the same block of $\Kr ( \pi )$
\end{tabular}  }
\right\} \ \Longrightarrow \ 
i_{k} = i_{l}.
\end{equation}
Under this assumption, we have to show that:
\begin{equation}
k_{\pi} ( \ a_{i_{n},i_{1}}^{(r_{1})}, a_{i_{1},i_{2}}^{(r_{2})}, 
\ldots , a_{i_{n-1},i_{n}}^{(r_{n})} \ ) \ =  \
[ \coef ( \ (r_{1},i_{1}), \ldots , (r_{n},i_{n}) \ ); \pi \ ] (f).
\end{equation}
Each of the two sides of (2.10) is a product of factors indexed by 
the blocks of $\pi$; we will prove (2.10) by showing that actually
for any given block $B$ of $\pi$, the factor corresponding to $B$ 
on the left-hand side of (2.10) is equal to the factor corresponding 
to $B$ on the right-hand side of (2.10).

So let us fix a block $B = \{ k_{1}<k_{2}< \cdots < k_{p} \}$ of $\pi$. 
The factor corresponding to $B$ on the left-hand side of (2.10) is:
\begin{equation}
k_{p} \Bigl(\ a_{ i_{k_{1}-1},i_{k_{1}} }^{ (r_{k_{1}}) },
a_{ i_{k_{2}-1},i_{k_{2}} }^{ (r_{k_{2}}) }, \ldots ,
a_{ i_{k_{p}-1},i_{k_{p}} }^{ (r_{k_{p}}) } \ \Bigr)
\end{equation}
(where recall that if $k_{1} = 1$, then we use $i_{n}$ for
``$i_{k_{1}-1}$''); the
factor corresponding to $B$ on the right-hand side of (2.10) is $[ \coef 
( \ (r_{k_{1}},i_{k_{1}}), \ldots , (r_{k_{p}},i_{k_{p}}) \ ) \ ] (f)$,
i.e:
\begin{equation}
k_{p} \Bigl( \ a_{ i_{k_{p}},i_{k_{1}} }^{ (r_{k_{1}}) },
a_{ i_{k_{1}},i_{k_{2}} }^{ (r_{k_{2}}) }, \ldots ,
a_{ i_{k_{p-1}},i_{k_{p}} }^{ (r_{k_{p}}) } \ \Bigr) .
\end{equation}

But now, let us notice that 
$k_{1}$ and $k_{2}-1$ belong to the same block of 
$\Kr ( \pi )$, and same for $k_{2}$ and $k_{3}-1, \ldots$ , same for
$k_{p}$ and $k_{1}-1$. This is easily seen by looking at the 
permutations associated to $\pi$ and $\Kr ( \pi )$: we have that
\[
\perm_{\pi} ( k_{1} ) = k_{2}, \ldots ,
\perm_{\pi} ( k_{p-1} ) = k_{p},
\perm_{\pi} ( k_{p} ) = k_{1},
\]
so from Eqn.(1.3) we get that:
\[
\perm_{\Kr ( \pi )} ( k_{2} -1 ) = k_{1}, \ldots ,
\perm_{\Kr ( \pi )} ( k_{p} -1 ) = k_{p-1},
\perm_{\Kr ( \pi )} ( k_{1} -1 ) = k_{p}.
\]
As a consequence of this remark and of the implication stated in (2.9), 
we see that the expressions appearing in (2.11) and (2.12)
are actually identical.

\vspace{6pt}

{\em Case 2:}
$[ \coef ( \ i_{1}, \ldots , i_{n} \ ); \Kr( \pi ) \ ] (G_{d}) =0.$

\vspace{6pt}

In this case we know that (2.9) does not hold, and we have to show that
the left-hand side of (2.8) is equal to 0.

It is immediate that, under the current assumption, we can find 
$1 \leq k,l \leq n$ such that:
\begin{equation}
\perm_{\Kr ( \pi )} (l) = k, \ \ \mbox{and} \ i_{k} \neq i_{l}.
\end{equation}
Indeed, if it were true that $i_{k} = i_{l}$ whenever 
$\perm_{\Kr ( \pi ) }(l) = k$, then by moving along the cycles of 
$\perm_{\Kr ( \pi ) }$ we would find that (2.9) holds.

By taking into account the relation between $\perm_{\pi}$ and
$\perm_{\Kr ( \pi )}$, we see that for $k,l$ as in (2.13) we also have 
that $\perm_{\pi} (k) = l+1$. Hence $k$ and $l+1$ belong to the same 
block $B$ of $\pi$; and moreover, if the block $B$ is written as
$B = \{ k_{1} < k_{2} < \cdots < k_{p} \}$, then there exists an index
$j$, $1 \leq j \leq p$ such that $k=k_{j}$ and $l+1 = k_{j+1}$ (with 
the convention that if $k=k_{p}$, then $l+1 = k_{1}$). But then the fact
that $i_{k} \neq i_{l}$ reads: $i_{k_{j}} \neq i_{k_{j+1}-1}$, which in 
turn implies that
\[
k_{p} \Bigl( \ a_{ i_{k_{1}-1},i_{k_{1}} }^{ (r_{k_{1}}) },
a_{ i_{k_{2}-1},i_{k_{2}} }^{ (r_{k_{2}}) }, \ldots ,
a_{ i_{k_{p}-1},i_{k_{p}} }^{ (r_{k_{p}}) } \ \Bigr) \ = \ 0
\]
(by the definition of R-cyclicity). Since the latter expression is the 
factor corresponding to $B$ in the product defining 
$k_{\pi} ( \ a_{i_{n},i_{1}}^{(r_{1})}, a_{i_{1},i_{2}}^{(r_{2})}, \ldots , 
a_{i_{n-1},i_{n}}^{(r_{n})} \ )$, we conclude that the left-hand side of 
(2.8) is indeed equal to 0. {\bf QED}

$\ $

{\bf Proof of Theorem 2.10.} Let $n$ be a positive integer, and consider 
some indices $1 \leq r_{1}, \ldots , r_{n} \leq s$,
$1 \leq i_{1}, \ldots , i_{n} \leq d$. By summing over $\pi \in NC(n)$
in the Equation (2.8) of Lemma 2.11, and by taking into account the 
properties of non-crossing cumulants and of boxed convolution, we get:
\begin{equation}
\varphi ( a_{i_{n},i_{1}}^{(r_{1})} a_{i_{1},i_{2}}^{(r_{2})} \cdots 
a_{i_{n-1},i_{n}}^{(r_{n})} ) \ = \ 
[ \coef ( \ (r_{1},i_{1}), \ldots , (r_{n},i_{n}) \ ) \ ]
( f \ \freestart \ G_{d} ).
\end{equation}
For every $1 \leq i \leq d$, let us denote by $P_{i} \in M_{d} ( \A )$
the matrix which has $I$ (the unit of $\A$) on the $(i,i)$-entry, and 
has all the other entries equal to 0. It is immediately verified that 
\[
\varphi_{d} ( A_{r_{1}}P_{i_{1}} \cdots A_{r_{n}}P_{i_{n}}  ) \ = \ 
\frac{1}{d} 
\varphi ( a_{i_{n},i_{1}}^{(r_{1})} a_{i_{1},i_{2}}^{(r_{2})} \cdots 
a_{i_{n-1},i_{n}}^{(r_{n})} ) .
\]
By combining this equation with (2.14), we get an equality of power 
series in $ds$ variables, which is stated as follows:
\begin{equation}
M_{A_{1}P_{1}, \ldots , A_{r}P_{i}, \ldots , A_{s}P_{d}} \ = \ 
\frac{1}{d} ( f \ \freestart \ G_{d} ).
\end{equation}
The Equation (2.6) is an immediate consequence of (2.15), since we have
for every $n \geq 1$ and $1 \leq r_{1}, \ldots , r_{n} \leq s$:
\[
\varphi_{d} ( A_{r_{1}} \cdots A_{r_{n}} ) \ = \ 
\sum_{i_{1}, \ldots , i_{n}=1}^{d}
\varphi_{d} ( A_{r_{1}}P_{i_{1}} \cdots A_{r_{n}}P_{i_{n}} )
\]
\[
= \ \sum_{i_{1}, \ldots , i_{n}=1}^{d}
[ \coef ( \ (r_{1},i_{1}), \ldots , (r_{n},i_{n}) \ ) \ ]
( M_{A_{1}P_{1}, \ldots , A_{s}P_{d}} )
\]
\[
= \ \frac{1}{d} \sum_{i_{1}, \ldots , i_{n}=1}^{d}
[ \coef ( \ (r_{1},i_{1}), \ldots , (r_{n},i_{n}) \ ) \ ]
( f \ \freestart \ G_{d} );
\]
the latter quantity is easily seen to be the coefficient of 
$z_{r_{1}} \cdots z_{r_{n}}$ in the series:
\[
\frac{1}{d} ( f \ \freestart \ G_{d} ) ( \ 
\underbrace{z_{1}, \ldots , z_{1},}_{d \ times} \ldots ,
\underbrace{z_{s}, \ldots , z_{s}}_{d \ times}  \ ),
\]
hence (2.6) follows.

On the other hand let us $\freestart$-convolve with 
$\mbox{M\"ob}_{d}$ on the right, on both sides of (2.15). On the 
left-hand side we get 
$M_{A_{1}P_{1}, \ldots , A_{s}P_{d}} \ \freestart \ \mbox{M\"ob}_{d}$,
which is equal to $ R_{A_{1}P_{1}, \ldots , A_{s}P_{d}}$ (see
Equation (1.20) in Section 1.5). On the right-hand side we get:
\[
\Bigl( \ \frac{1}{d} ( f \ \freestart \ G_{d} ) \ \Bigr) \ 
\freestart \ \mbox{M\"ob}_{d} 
\ = \ \frac{1}{d} \Bigl( \ f \ \freestart \ G_{d} \ \freestart \ 
( d \mbox{M\"ob}_{d} \circ D_{1/d}) \ \Bigr)
\ \ \mbox{ (by Eqn.(1.23) )}
\]
\[
= \ \frac{1}{d} ( f \ \freestart \ H_{d} )
\ \ \mbox{ (by the definition of $H_{d}$ in Section 1.7).}
\]
So we obtain the equation:
\begin{equation}
R_{A_{1}P_{1}, \ldots , A_{r}P_{i}, \ldots , A_{s}P_{d}} \ = \ 
\frac{1}{d} ( f \ \freestart \ H_{d} ),
\end{equation}
out of which (2.7) is obtained in the same way as (2.6) was obtained 
from (2.15). {\bf QED}

$\ $

{\bf 2.12 Remark.} The proof of Theorem 2.10 obtains the Equations
(2.15) and (2.16), stronger than what was originally stated, and which
show better the significance of the series $f \ \freestart \ G_{d}$ 
and $f \ \freestart \ H_{d}$.

$\ $

$\ $

\setcounter{section}{3}
{\large\bf 3. Applications of Theorem 2.10.}

$\ $

We will concentrate on applications to a family $A_{1}, \ldots , A_{s}$ 
\setcounter{equation}{0}
of selfadjoint $d \times d$ matrices over a $*$-probability space $\ncps$.
By keeping in mind the motivating example from \cite{V2}, it is of 
particular interest to put into evidence situations where the family 
$A_{1}, \ldots , A_{s}$ is free in $(M_{d} ( \A ), \varphi_{d} )$, and 
where the individual R-transform of each of $A_{1}, \ldots , A_{s}$ is 
determined explicitly. It seems that some important situations of this 
kind appear as a consequence of a ``partial summation condition'',
described in the next proposition.

$\ $

{\bf 3.1 Proposition.}  Let $\ncps$ be a $*$-probability space, let $d,s$
be positive integers, and let $A_{1} = [ a_{i,j}^{(1)} ]_{i,j=1}^{d},
\ldots , A_{s} = [ a_{i,j}^{(s)} ]_{i,j=1}^{d}$ form an R-cyclic
family of selfadjoint matrices in $M_{d} ( \A )$. We denote the 
determining series of $A_{1}, \ldots , A_{s}$ by $f$. Suppose that 
for every $n \geq 1$ and every $1 \leq r_{1}, \ldots, r_{n} \leq s$,
$1 \leq i_{1}, \ldots , i_{n} \leq d$, the sum:
\begin{equation}
\sum_{i_{1}, \ldots , i_{n-1}=1}^{d} \ 
[ \coef ( \ (r_{1},i_{1}), \ldots , (r_{n-1},i_{n-1}), (r_{n},i_{n}) \ ) ]
(f) \ =: \lambda_{r_{1}, \ldots , r_{n}}
\end{equation}
does not depend on $i_{n}$ (even though the sum is only over 
$i_{1}, \ldots , i_{n-1}$). Then:
\begin{equation}
R_{A_{1}, \ldots , A_{s}} ( z_{1}, \ldots , z_{s} ) \ = \
\sum_{n=1}^{\infty} \ \sum_{r_{1}, \ldots , r_{n}=1}^{s} \ 
\lambda_{r_{1}, \ldots , r_{n}} z_{r_{1}} \cdots z_{r_{n}}.
\end{equation}

$\ $

{\bf Proof.} The Equation (3.2) is equivalent to the fact that for every 
$n \geq 1$ and every $1 \leq r_{1}, \ldots , r_{n} \leq s$ we have:
\begin{equation}
k_{n} ( A_{r_{1}}, \ldots , A_{r_{n}} ) \ = \
\lambda_{r_{1}, \ldots , r_{n}}.
\end{equation}
We fix $n$ and $r_{1}, \ldots , r_{n}$ about which we show that (3.3)
is true. The case when $n=1$ is trivial, so we will assume that $n \geq 2$.

The Equation (2.7) of Theorem 2.10 gives us the formula:
\[
k_{n} ( A_{r_{1}}, \ldots , A_{r_{n}} ) \ = 
\]
\[
\frac{1}{d} \ \sum_{i_{1}, \ldots , i_{n}=1}^{d} \ \sum_{\pi \in NC(n)} \
[ \coef ( \ (r_{1},i_{1}), \ldots , (r_{n},i_{n}) \ ); \pi ] (f) \cdot
[ \coef ( i_{1}, \ldots , i_{n}) ; \Kr ( \pi ) ] (H_{d}).
\]
We will write this in the form:
\begin{equation}
k_{n} ( A_{r_{1}}, \ldots , A_{r_{n}} ) \ = \ \sum_{\pi \in NC(n)}
T_{\pi},
\end{equation}
where for every $\pi \in NC(n)$ we set:
\begin{equation}
T_{\pi} \ := \ \frac{1}{d} \ \sum_{i_{1}, \ldots , i_{n}=1}^{d} \ 
[ \coef ( \ (r_{1},i_{1}), \ldots , (r_{n},i_{n}) \ ); \pi ] (f) \cdot
[ \coef ( i_{1}, \ldots , i_{n}) ; \Kr ( \pi ) ] (H_{d}).
\end{equation}

We first consider the quantity $T_{\pi}$ defined in (3.5) in the special
case when $\pi = 1_{n}$, the partition of $\{ 1, \ldots , n \}$ which 
has only one block. In this case $\Kr ( \pi )$ is the partition into $n$
blocks of one element; since all the coefficients of degree 1 of $H_{d}$
are equal to 1, it follows that
\[
[ \coef ( i_{1}, \ldots , i_{n}) ; \Kr ( 1_{n} ) ] (H_{d}) \ = \ 1, \ \ 
\forall \ i_{1}, \ldots , i_{n} \in \{ 1, \ldots , d \} .
\]
We hence get:
\[
T_{1_{n}} \ = \ \frac{1}{d} \ \sum_{i_{1}, \ldots , i_{n} =1}^{d} \
[ \coef ( \ (r_{1},i_{1}), \ldots , (r_{n},i_{n}) \ ) ] (f).
\]
The partial summation property of the series $f$ (given in Eqn.(3.1))
implies that the latter sum is equal to $\lambda_{r_{1}, \ldots , r_{n}}$.
Thus, in view of (3.4), the proof will be over if we can show that 
$T_{\pi} = 0$ for every $\pi \neq 1_{n}$ in $NC(n)$. 

So for the remaining of the proof we fix a partition $\pi \neq 1_{n}$
in $NC(n)$. Moreover, we will also fix a block $B_{o}$ of $\pi$ which is
an interval, $B_{o} = [ p, q ] \cap \Z$ with $1 \leq p \leq q \leq n$
(every non-crossing partition has such a block). The considerations below,
leading to the conclusion that $T_{\pi} = 0$, will be made by looking at
the case when $B_{o}$ has more than one element; the case when 
$|B_{o}|=1$ (which is similar, and easier) is left as an exercise to the
reader. We denote by ``$Rest$'' the set of blocks of $\pi$ which are
different from $B_{o}$.

Let us now look at at the Kreweras complement $\Kr ( \pi )$. It is 
immediate that $\{ p \}$, $\{ p+1 \} , \ldots , \{ q-1 \}$ are 
one-element blocks of $\Kr ( \pi )$. We denote by $B_{o}'$ the block of 
$\Kr ( \pi )$ which contains $q$; observe that $B_{o}'$ has more than
one element -- indeed, it is clear that $p-1$ also belongs to 
$B_{o}'$ (where if $p=1$, then ``$p-1$'' means ``$n$''; even in this case 
we have that $p-1 \neq q$, since it was assumed that $\pi \neq 1_{n}$).
Moreover, let us denote by $Rest'$ the set of blocks of $\Kr ( \pi )$
(if any) which remain after $\{ p \}$, $\{ p+1 \} , \ldots , \{ q-1 \}$
and $B_{o}'$ are deleted.

For any $i_{1}, \ldots , i_{n} \in \{ 1, \ldots , d \}$ we have:
\[
[ \coef ( \ (r_{1},i_{1}), \ldots , (r_{n},i_{n}) \ ); \pi ] (f) \cdot
[ \coef ( i_{1}, \ldots , i_{n}) ; \Kr ( \pi ) ] (H_{d}) \ = 
\]
\begin{equation}
[ \coef ( \ (r_{p},i_{p}), \ldots , (r_{q},i_{q}) \ ) ] (f) \cdot
[ \coef ( i_{1}, \ldots , i_{n}) | B_{o}' ] (H_{d}) \cdot
\end{equation}
\[
\cdot \prod_{B \in Rest} 
[ \coef ( \ (r_{1},i_{1}), \ldots , (r_{n},i_{n}) \ ) | B ] (f) \cdot
\prod_{B' \in Rest'} 
[ \coef ( i_{1}, \ldots , i_{n}) | B' ] (H_{d}) 
\]
(we took into account that the factors
$[ \coef (i_{p}) ] (H_{d}), \ldots , [ \coef (i_{q-1}) ] (H_{d})$,
which should also appear on the right-hand side of (3.6), are all equal 
to 1). The indices $i_{p}, \ldots , i_{q-1}$ appear only in the factor 
``$[ \coef ( \ (r_{p},i_{p}), \ldots , (r_{q},i_{q}) \ ) ] (f)$'' of 
(3.6). Thus, if in (3.6) we sum over $i_{p}, \ldots , i_{q-1}$, and 
make use of the partial summation property from (3.1), then we get:
\begin{equation}
\lambda_{r_{p}, \ldots , r_{q}} \cdot
[ \coef ( i_{1}, \ldots , i_{n}) | B_{o}' ] (H_{d}) \cdot
\end{equation}
\[
\cdot \prod_{B \in Rest} 
[ \coef ( \ (r_{1},i_{1}), \ldots , (r_{n},i_{n}) \ ) | B ] (f) \cdot
\prod_{B' \in Rest'} 
[ \coef ( i_{1}, \ldots , i_{n}) | B' ] (H_{d})
\]
(expression depending on some arbitrary indices 
$i_{1}, \ldots , i_{p-1}, i_{q}, \ldots , i_{n}$, chosen from
$\{ 1, \ldots , d \}$).

Next, in (3.7) we sum over the index $i_{q}$. The only factor in (3.7)
which involves $i_{q}$ is
``$[ \coef ( i_{1}, \ldots , i_{n}) | B_{o}' ] (H_{d})$'', so as a result
of this new summation we get:
\[
\lambda_{r_{p}, \ldots , r_{q}} \cdot \Bigl\{ \ \sum_{i_{q} = 1}^{d}
[ \coef ( i_{1}, \ldots , i_{n}) | B_{o}' ] (H_{d}) \ \Bigr\} \cdot
\]
\[
\cdot \prod_{B \in Rest} 
[ \coef ( \ (r_{1},i_{1}), \ldots , (r_{n},i_{n}) \ ) | B ] (f) \cdot
\prod_{B' \in Rest'} 
[ \coef ( i_{1}, \ldots , i_{n}) | B' ] (H_{d}).
\]
But, as an immediate consequence of the remark concluding Section 1.7, 
we have that $\sum_{i_{q} = 1}^{d}
[ \coef ( i_{1}, \ldots , i_{n}) | B_{o}' ] (H_{d}) = 0$.

The conclusion that we draw from the preceding three paragraphs is the
following: for any choice of the indices
$i_{1}, \ldots , i_{p-1}, i_{q+1}, \ldots , i_{n} \in \{ 1, \ldots , d \}$,
we have that
\[
\sum_{i_{p}, \ldots , i_{q}=1}^{d} \ 
[ \coef ( \ (r_{1},i_{1}), \ldots , (r_{n},i_{n}) \ ); \pi ] (f) \cdot
[ \coef ( i_{1}, \ldots , i_{n}) ; \Kr ( \pi ) ] (H_{d}) \ = \ 0.
\]
It only remains  that we sum over
$i_{1}, \ldots , i_{p-1}, i_{q+1}, \ldots , i_{n}$ in the latter equation,
to obtain the desired fact that $T_{\pi} = 0$. {\bf QED}

$\ $

{\bf 3.2 Corollary.}  Let $\ncps$ be a $*$-probability space, let $d,s$
be positive integers, and let $A_{1} = [ a_{i,j}^{(1)} ]_{i,j=1}^{d},
\ldots , A_{s} = [ a_{i,j}^{(s)} ]_{i,j=1}^{d}$ form an R-cyclic
family of selfadjoint matrices in $M_{d} ( \A )$. Suppose that the $s$
families of entries $\{ a_{i,j}^{(r)} \ | \ 1 \leq i,j \leq d \}$, 
with $1 \leq r \leq s$, are free in $\ncps$. Moreover, for every 
$1 \leq r \leq s$ let $f_{r} \in \Theta_{d}$ be the determining series 
of $A_{r}$. We assume that for every $n \geq 1$ and for every 
$1 \leq r \leq s$, $1 \leq i \leq d$, the sum:
\begin{equation}
\sum_{i_{1}, \ldots , i_{n-1}=1}^{d} \
[ \coef (i_{1}, \ldots , i_{n-1},i) ] (f_{r}) \ =: \ \lambda_{n}^{(r)}
\end{equation}
does not depend on the choice of $i$ (but only on $n$ and $r$). Then the 
matrices $A_{1}, \ldots , A_{s}$ are free in $(M_{d} ( \A ), \varphi_{d})$,
and have R-transforms
\begin{equation}
R_{A_{r}} (z) \ = \ \sum_{n=1}^{\infty} \lambda_{n}^{(r)} z^{n}, \ \ 
1 \leq r \leq s.
\end{equation}

$\ $

{\bf Proof.} Let $f$ denote the determining series of the whole R-cyclic 
family $A_{1}, \ldots , A_{s}$. The condition of freeness between the
families of entries of $A_{1}, \ldots , A_{s}$ implies the formula:
\[
f( z_{1,1}, \ldots , z_{r,i}, \ldots , z_{s,d} ) \ = \ 
\sum_{r=1}^{s} f_{r} ( z_{r,1}, \ldots , z_{r,i}, \ldots , z_{r,d} ),
\]
where $f_{r}$ is (as in the statement of the corollary) the 
determining series for just the R-cyclic matrix $A_{r}$. It is immediate
that $f$ satisfies the partial summation condition described in 
Equation (3.1) of Proposition 3.1, where we set:
\[
\lambda_{r_{1}, \ldots , r_{n}} \ = \ \left\{  
\begin{array}{cl}
\lambda_{n}^{(r)}  &  \mbox{if }  r_{1} = \cdots = r_{n} = r  \\
0                  &  \mbox{otherwise.}
\end{array}  \right.
\]
Thus the Proposition 3.1 can be applied, and gives us:
\[
R_{A_{1}, \ldots , A_{s}} (z_{1}, \ldots , z_{s}) \ = \ 
\sum_{r=1}^{s} \ \sum_{n=1}^{\infty} \lambda_{n}^{(r)} z_{r}^{n},
\]
which (by virtue of Equation (1.15) in Section 1) is equivalent to saying 
that $A_{1}, \ldots , A_{s}$ are free and and have the indicated individual 
R-transforms. {\bf QED}

$\ $

The Corollary 3.2 can be in turn particularized to the situation of a family
of matrices with free R-diagonal entries (on the line of Example 2.6). The 
precise spelling of this particular case goes as follows.

$\ $

{\bf 3.3 Corollary.}  Let $\ncps$ be a $*$-probability space, let $d,s$
be positive integers, and suppose that the elements 
$\{ a_{i,j}^{(r)} \ | \ 1 \leq i,j \leq d, \ 1 \leq r \leq s \}$ of $\A$
have the following properties:

(i) For every $1 \leq i \leq d$ and $1 \leq r \leq s$, the element 
$a_{i,i}^{(r)}$ is selfadjoint. We denote the R-transform of 
$a_{i,i}^{(r)}$ as $\sum_{n=1}^{\infty} \alpha_{i,i;n}^{(r)} z^{n}$.

(ii) For every $1 \leq i,j \leq d$ such that $i \neq j$, and for every 
$1 \leq r \leq s$, the element $a_{i,j}^{(r)}$ is R-diagonal and has 
$\Bigl( a_{i,j}^{(r)} \Bigr)^{*} = a_{j,i}^{(r)}$. We denote the 
determining series of $a_{i,j}^{(r)}$ as 
$\sum_{n=1}^{\infty} \alpha_{i,j;2n}^{(r)} z^{n}$;
we also set $\alpha_{i,j;2n-1}^{(r)} := 0$, $\forall \ n \geq 1$.

(iii) The $sd(d+1)/2$ families: 
$\{ a_{i,i}^{(r)} \}$ for $1 \leq i \leq d$, $1 \leq r \leq s$, together 
with $\{ a_{i,j}^{(r)}, a_{j,i}^{(r)} \}$ for $1 \leq i<j \leq d$, 
$1 \leq r \leq s$ are free in $\ncps$.

Suppose moreover that for every $n \geq 1$ and every $1 \leq r \leq s$, 
$1 \leq i \leq d$, the sum:
\begin{equation}
\sum_{j=1}^{d} \ \alpha_{i,j;n}^{(r)} \ =: \ \lambda_{n}^{(r)}
\end{equation}
does not actually depend on $i$. Then the 
matrices $A_{1} = [ a_{i,j}^{(1)} ]_{i,j=1}^{d}, \ldots ,
A_{s} = [ a_{i,j}^{(s)} ]_{i,j=1}^{d}$ are free in 
$(M_{d} ( \A ), \varphi_{d})$, and have R-transforms
\[
R_{A_{r}} (z) \ = \ \sum_{n=1}^{\infty} \lambda_{n}^{(r)} z^{n}, \ \ 
1 \leq r \leq s.
\]

$\ $

{\bf 3.4 Remark.} The summation conditions (3.10) become extremely simple
when the elements $a_{i,i}^{(r)}$ are semicircular, and the elements 
$a_{i,j}^{(r)}$ with $i \neq j$ are circular. Indeed, in this case 
we have that $\alpha_{i,j;n}^{(r)} = 0$ whenever $n \neq 2$, and that
$\alpha_{i,j;2}^{(r)}$ is one quarter of the squared radius of the 
circular/semicircular element $a_{i,j}^{(r)}$ (compare to Remark 2.7.2).
Thus in this case if we denote the radius of $a_{i,j}^{(r)}$ by 
$\gamma_{i,j}^{(r)}$, then (3.10) amounts to asking that for every 
$1 \leq r \leq s$ the matrix of squared radii
$[ \gamma_{i,j}^{(r)} ]_{i,j=1}^{d}$ has constant sums along its columns:
\[
\sum_{j=1}^{d} \Bigl( \gamma_{1,j}^{(r)} \Bigr)^{2} \ = \ \cdots \ = \
\sum_{j=1}^{d} \Bigl( \gamma_{d,j}^{(r)} \Bigr)^{2} \ =: \ 
\gamma_{r}^{2}.
\]
The conclusion of Corollary 3.3 becomes that the matrices
$A_{1} = [ a_{i,j}^{(1)} ]_{i,j=1}^{d}, \ldots ,
A_{s} = [ a_{i,j}^{(s)} ]_{i,j=1}^{d}$ are free, and that $A_{r}$ is 
semicircular of radius $\gamma_{r}$, for $1 \leq r \leq s$. This particular
case of Corollary 3.3 is very close to Proposition 2.9 of
\cite{V2}, and can also be obtained by the methods used there 
(approximations in distribution by large Gaussian random matrices).

$\ $

Another particularization of Proposition 3.1 covers a situation when
the matrices $A_{1}, \ldots$, 
$A_{s}$ are not free, and which is 
motivated by results about free compressions (see Sections 8 and 17 of
\cite{NS3}; the case of only one matrix appeared in Example 2.4 above).

$\ $

{\bf 3.5 Corollary.} Let $\ncps$ be a $*$-probability space, let $d,s$
be positive integers, and let $A_{1} = [ a_{i,j}^{(1)} ]_{i,j=1}^{d},
\ldots , A_{s} = [ a_{i,j}^{(s)} ]_{i,j=1}^{d}$ form an R-cyclic
family of selfadjoint matrices in $M_{d} ( \A )$. Suppose that the 
cyclic cumulants of the entries of these matrices depend only on the
superscript indices:
\begin{equation}
k_{n}( a_{i_{n},i_{1}}^{(r_{1})}, a_{i_{1},i_{2}}^{(r_{2})}, 
\ldots , a_{i_{n-1},i_{n}}^{(r_{n})} ) \ =: \
\alpha_{ r_{1}, \ldots , r_{n} },
\end{equation}
for every $n \geq 1$ and every $1 \leq r_{1}, \ldots , r_{n} \leq s$,
$1 \leq i_{1}, \ldots , i_{n} \leq d$. Then:
\begin{equation}
R_{A_{1}, \ldots , A_{s}} (z_{1}, \ldots , z_{s}) \ = \ 
\sum_{n=1}^{\infty} \ \sum_{r_{1}, \ldots , r_{n}=1}^{s} d^{n-1}
\alpha_{ r_{1}, \ldots , r_{n} } z_{r_{1}} \cdots z_{r_{n}}.
\end{equation}

$\ $

{\bf Proof.} If $f$ denotes the determining series of 
$A_{1}, \ldots , A_{s}$, then the coefficients of $f$ are:
\[
[ \coef ( \ (r_{1},i_{1}), \ldots , (r_{n},i_{n}) \ ) ] (f) \ =: 
\alpha_{r_{1}, \ldots , r_{n}},
\]
$\forall \ n \geq 1$, $\forall \ 1 \leq r_{1}, \ldots , r_{n} \leq s$,
$1 \leq i_{1}, \ldots , i_{n} \leq d$.
It is obvious that the partial summation condition of Equation (3.1)
holds, where $\lambda_{r_{1}, \ldots , r_{n}} := d^{n-1}
\alpha_{r_{1}, \ldots , r_{n}}$. {\bf QED}

$\ $

$\ $

\setcounter{section}{4}
{\large\bf 4. Algebras generated by R-cyclic families}

$\ $

{\bf 4.1 Remark.} Let $\ncps$ be a non-commutative probability space,
\setcounter{equation}{0}
let $d$ be a positive integer, and let $A_{1}, \ldots , A_{s}$ be an 
R-cyclic family of matrices in $M_{d} ( \A )$. Directly from the 
definition of R-cyclicity, and by using some basic properties of the
non-crossing cumulants, it is easy to observe several ``operations'' that 
can be performed on the family $A_{1}, \ldots , A_{s}$ without affecting
its R-cyclicity. For instance, it is trivial that re-ordering the $s$ 
matrices does not affect R-cyclicity, and same about the operation of 
deleting one of the matrices from the family. Another operation which
clearly does not affect the R-cyclicity of $A_{1}, \ldots , A_{s}$ consists
in arbitrarily re-scaling the entries of the matrices (multiply the 
$(i,j)$-entry of $A_{r}$ by some constant $\lambda_{i,j}^{(r)}$, for every
$1 \leq i,j \leq d$, $1 \leq r \leq s$). Let us also observe that:

\vspace{6pt}

(a) If we enlarge $A_{1}, \ldots , A_{s}$ with a matrix 
$A \in \mbox{span} \{ A_{1}, \ldots , A_{s} \}$, then the enlarged
family $A_{1}, \ldots , A_{s}, A$ is still R-cyclic. This is a direct
consequence of the multilinearity of the cumulant functionals 
$k_{n} : \A^{n} \rightarrow \C$, $n \geq 1$.

\vspace{6pt}

(b) If we enlarge $A_{1}, \ldots , A_{s}$ with a scalar diagonal matrix 
$D$ (which has the diagonal entries of the form $\lambda_{i} I$,
$1 \leq i \leq d$, and the off-diagonal entries equal to 0), then the
enlarged family $A_{1}, \ldots , A_{s}, D$ is still R-cyclic. This is
a consequence of the fact that a non-crossing cumulant of $n \geq 2$
variables is 0 if at least one of its entries is in $\C I$ (same kind of
argument as in the last phrase of Section 1).

$\ $

In connection to (b) of Remark 4.1, note that we could not use there a 
scalar matrix which is not diagonal -- indeed, the R-cyclicity condition 
asks in particular that every off-diagonal entry of every matrix 
in the family lies in the kernel of the state $\varphi$.

Now, in the framework of the same R-cyclic family $A_{1}, \ldots , A_{s}$
as above, where we assume that $s \geq 2$, let us also observe that:

$\ $

{\bf 4.2 Lemma.} If $A_{s+1} := A_{1}A_{2}$, then the enlarged family
$A_{1}, \ldots , A_{s}, A_{s+1}$ is still R-cyclic.

$\ $

{\bf Proof.} We will use a formula for free cumulants with products as
entries, as developed in \cite{KS}. In fact we can set the proof by 
induction, in such a way that we only use a particular case of
this formula, which had already appeared in \cite{S2}. The particular
case in question says that for any $1 \leq m < n$ and any
$x_{1}, \ldots , x_{n}$ in $\A$ we have:
\[
k_{n-1} ( x_{1}, \ldots , x_{m-1}, x_{m}x_{m+1}, 
x_{m+2}, \ldots , x_{n} )
\]
\[
\ = \ k_{n} ( x_{1}, \ldots , x_{n} ) + 
k_{m} ( x_{1}, \ldots , x_{m} ) \cdot
k_{n-m} ( x_{m+1}, \ldots , x_{n} )
\]
\begin{equation}
+ \sum_{j=2}^{m} 
k_{m-j+1} ( x_{j}, \ldots , x_{m} ) \cdot
k_{n-m+j-1} ( x_{1}, \ldots , x_{j-1}, x_{m+1}, \ldots , x_{n} )
\end{equation}
\[
+ \sum_{j=m+1}^{n-1} 
k_{j-m} ( x_{m+1}, \ldots , x_{j} ) \cdot
k_{n-j+m} ( x_{1}, \ldots , x_{m}, x_{j+1}, \ldots , x_{n} ) .
\]

Now, let us return to the matrices $A_{1}, \ldots , A_{s+1}$ appearing
in the statement of the lemma. For $1 \leq r \leq s+1$ and
$1 \leq i,j \leq d$ we denote by $a_{i,j}^{(r)}$ the $(i,j)$-entry of 
$A_{r}$. The hypothesis that $A_{s+1} = A_{1}A_{2}$ thus says that
\begin{equation}
a_{i,j}^{(s+1)} \ = \ \sum_{k=1}^{d}
a_{i,k}^{(1)} a_{k,j}^{(2)} , \ \ \forall \ 1 \leq i,j \leq d.
\end{equation}

We will prove by induction on $l \geq 0$ the following statement:
\[
\mbox{St} (l) \ \ 
\left\{  \begin{array}{l}
\mbox{For every $n \geq 1$, $r_{1}, \ldots , r_{n} \in 
\{ 1, \ldots , s+1 \}$ and $i_{1},j_{1}, \ldots , i_{n},j_{n}
\in \{ 1, \ldots , d \}$ }                                       \\ 

\mbox{such that $\{ m \ | \ 1 \leq m \leq n, \ r_{m}
= s+1 \}$ has $l$ elements}                                        \\

\mbox{and for which it is not true that $j_{1}=i_{2}, \ldots , 
j_{n-1}=i_{n}, j_{n}=i_{1}$ }                                     \\

\mbox{we have that $k_{n} ( a_{i_{1},j_{1}}^{(r_{1})}, \ldots ,
a_{i_{n},j_{n}}^{(r_{n})} ) = 0$. }
\end{array}  \right.
\]
If $l=0$, the statement $\mbox{St}(l)$ amounts precisely to the hypothesis 
that the family $A_{1}, \ldots , A_{s}$ is R-cyclic. For the rest of the 
proof we fix an $l \geq 1$, for which we assume that the statements
$\mbox{St}(0), \ldots , \mbox{St}(l-1)$ are true, and for which we prove 
that the statement $\mbox{St}(l)$ is also true.

Consider $n \geq 1$, $r_{1}, \ldots , r_{n} \in \{ 1, \ldots , s+1 \}$
and $i_{1},j_{1}, \ldots , i_{n},j_{n} \in \{ 1, \ldots , d \}$ such 
that $\{ m \ | \ 1 \leq m \leq n, \ r_{m} = s+1 \}$ has $l$ elements,
and for which it is not true that
$j_{1}=i_{2}, \ldots , j_{n-1}=i_{n}, j_{n}=i_{1}$. Moreover, let us fix
an index $m$, $1 \leq m \leq n$, such that $r_{m} = s+1$. Our goal is 
to show that $k_{n} ( a_{i_{1},j_{1}}^{(r_{1})}, \ldots , 
a_{i_{n},j_{n}}^{(r_{n})} ) = 0$, but in view of (4.2) and of the 
multilinearity of $k_{n}$ it suffices to verify that: 
\begin{equation}
k_{n} ( a_{i_{1},j_{1}}^{(r_{1})}, \ldots , 
a_{i_{m-1},j_{m-1}}^{(r_{m-1})}, 
a_{i_{m},k}^{(1)} a_{k,j_{m}}^{(2)},
a_{i_{m+1},j_{m+1}}^{(r_{m+1})}, \ldots , 
a_{i_{n},j_{n}}^{(r_{n})} ) = 0, \ \ \forall \ 1 \leq k \leq d.
\end{equation}

Finally, let us also fix an index $k \in \{ 1, \ldots , d \}$ about which 
we will show that (4.3) holds. This in fact will be an immediate application
of the formula (4.1). Indeed, let us pick an index 
$p \in \{ 1, \ldots , n \}$ such that $j_{p} \neq i_{p+1}$; for the sake
of clarity of the presentation we will assume that we know the relative
position of $p$ and $m$ -- say for instance that $p<m-1$ (all the cases are 
treated similarly). We apply the formula (4.1) to the cumulant (4.3), and 
obtain a sum of $n+1$ terms $T_{1}, T_{2}, \ldots , T_{n+1}$ where each 
of these terms is either a cumulant or a product of two cumulants:
\begin{equation}
k_{n} ( a_{i_{1},j_{1}}^{(r_{1})}, \ldots , 
a_{i_{p},j_{p}}^{(r_{p})}, 
a_{i_{p+1},j_{p+1}}^{(r_{p+1})}, \ldots ,
a_{i_{m},k}^{(1)} a_{k,j_{m}}^{(2)}, \ldots , 
a_{i_{n},j_{n}}^{(r_{n})} ) 
\end{equation}
\[
= \ T_{1} + T_{2} + \cdots + T_{n+1} .
\]
The list of superscript indices on the left-hand side of (4.4) is 
$r_{1}, \ldots , r_{m-1}, 1,2, r_{m+1}, \ldots , r_{n}$, containing
$l-1$ occurrences of $s+1$. So the induction hypothesis will apply and
will give us that $T_{1} = \cdots = T_{n+1} = 0$ on the right-hand side 
of (4.4), provided that we make sure that each of $T_{1}, \ldots , T_{n+1}$
still violates the cyclicity condition of the subscript indices. The
violation of cyclicity for subscript indices is trivial for all of 
$T_{1}, \ldots , T_{n+1}$ with one exception, because in general the 
neighboring indices $j_{p} \neq i_{p+1}$ will not be separated. The 
exception is for the term:
\[
k_{m-p} ( a_{i_{p+1},j_{p+1}}^{(r_{p+1})}, \ldots ,
a_{i_{m},k}^{(1)} ) \cdot
k_{n+1-m+p} ( a_{i_{1},j_{1}}^{(r_{1})}, \ldots ,
a_{i_{p},j_{p}}^{(r_{p})}, a_{k,j_{m}}^{(2)}, \ldots ,
a_{i_{n},j_{n}}^{(r_{n})} ) ;
\]
but here the cyclicity condition of the subscript indices is still 
violated, since we must have that either $k \neq i_{p+1}$ or that 
$j_{p} \neq k$. {\bf QED}

$\ $

By combining the various ``R-cyclicity preserving operations'' which
were observed in the Remark 4.1 and Lemma 4.2, we arrive to the 
following statement (which in some sense collects these observations
together):

$\ $

{\bf 4.3 Theorem.}  Let $\ncps$ be a non-commutative probability space,
let $d$ be a positive integer, and let $A_{1}, \ldots , A_{s}$ be an 
R-cyclic family of matrices in $M_{d} ( \A )$. We denote by $\D$ the 
algebra of scalar diagonal matrices in $M_{d} ( \A )$, and by 
${\cal C}$ the subalgebra of $M_{d} ( \A )$ which is generated by
$\{ A_{1}, \ldots , A_{s} \} \cup \D$.
Then every finite family of matrices from ${\cal C}$ is R-cyclic.

$\ $

This theorem will be put into a better perspective by the result
in Section 8.

$\ $

$\ $

\setcounter{section}{5}
{\large\bf 5. Review of operator-valued cumulants}

$\ $

Let $\B$ be a unital algebra over $\C$. By a {\em $\B$-probability space}
\setcounter{equation}{0}
we understand a pair $( \M , E )$, where:

-- $\M$ is an algebra containing $\B$ as a unital subalgebra (by which 
we mean that $\B$ is identified as a unital subalgebra of $\M$, in some
well-determined way); 

-- $E: \M \rightarrow \B$ is a linear map with the properties that
$E(b) = b$ for every $b \in \B$, and $E(b_{1}xb_{2}) = b_{1}E(x)b_{2}$ for
every $b_{1}, b_{2} \in \B$, $x \in \M$.

If $( \M , E)$ is a $\B$-probability space and if 
$x_{1}, \ldots , x_{s} \in \M$, then the expressions of the form:
\[
E( b_{0}x_{r_{1}}b_{1} \cdots x_{r_{n}}b_{n} ), \mbox{ with }
n \geq 1, \ 1 \leq r_{1}, \ldots , r_{n} \leq s, \
b_{0}, b_{1}, \ldots , b_{n} \in \B
\]
are called {\em joint $\B$-moments} of the family $x_{1}, \ldots , x_{s}$.
Moreover, if $( \widetilde{\M}, \widetilde{E} )$ also is a $\B$-probability
space and if 
$\widetilde{x}_{1}, \ldots , \widetilde{x}_{s} \in \widetilde{\M}$, 
we will say that the families $x_{1}, \ldots , x_{s}$ and 
$\widetilde{x_{1}}, \ldots , \widetilde{x_{s}}$ have
{\em identical $\B$-distributions} if
\begin{equation}
E( b_{0}x_{r_{1}}b_{1} \cdots x_{r_{n}}b_{n} ) \ = \ 
\widetilde{E} ( b_{0}\widetilde{x}_{r_{1}}b_{1} \cdots 
\widetilde{x}_{r_{n}}b_{n} ) 
\end{equation}
for every $n \geq 1, \ 1 \leq r_{1}, \ldots , r_{n} \leq s,$ and
$b_{0}, b_{1}, \ldots , b_{n} \in \B$.

While the joint $\B$-moments generalize the joint moments appearing in
Eqn.(1.1) of Section 1, it is in general inconvenient to introduce a 
concept of ``$\B$-moment series'' analogous to the one defined by 
Eqn.(1.2). Similarly, rather than introducing $\B$-valued R-transforms, 
it is more convenient to just consider the $\B$-valued generalization for
the concept of non-crossing cumulant. Following the development of 
\cite{S2}, this can be done as described in Proposition 5.2 below.

$\ $

{\bf 5.1 Notations.} Let $\pi , \rho$ be partitions in $NC(p)$ and 
$NC(q)$ respectively, where $p,q \geq 1$. Let $k$ be in 
$\{ 0,1, \ldots , q \}$. By $\ins ( \pi \mapsto \rho ; k )$ we will
denote the non-crossing partition in $NC(p+q)$ which is obtained by
``inserting $\pi$ between the elements $k$ and $k+1$ of $\rho$''.
Formally this means that the set $\{ k+1, \ldots , k+p \}$ is a 
union of blocks of $\ins ( \pi \mapsto \rho ; k )$, and that:

(i) the restriction of $\ins ( \pi \mapsto \rho ; k )$ to 
$\{ k+1, \ldots , k+p \}$ is naturally identified to $\pi$;

(ii) the restriction of $\ins ( \pi \mapsto \rho ; k )$ to 
$\{ 1, 2, \ldots , p+q \} \setminus \{ k+1, \ldots , k+p \}$ is naturally 
identified to $\rho$.

$\ $

For example, if $\pi = \{ \ \{ 1 \} , \{ 2,3 \} \ \} \in NC(3)$ and 
$\rho = \{ \ \{ 1,2  \} \ \} \in NC(2)$, then: 
$\ins ( \pi \mapsto \rho ; 0)$ = 
$\{ \ \{ 1 \} ,  \{ 2,3 \} , \{ 4,5  \} \ \}$;
$\ins ( \pi \mapsto \rho ; 1)$ = $\{ \ \{ 1,5 \} , \{ 2 \} , \{ 3,4 \} \ \}$;
$\ins ( \pi \mapsto \rho ; 2)$ = $\{ \ \{ 1,2 \} , \{ 3 \} , \{ 4,5 \} \ \}$.

$\ $

{\bf 5.2 Proposition} (see \cite{S2}, Section 3.2).
Let $( \M , E)$ be a $\B$-probability space. There exists a family of 
functionals $\{ k_{\pi}^{( \B )} \ | \ \pi \in \cup_{n=1}^{\infty} NC(n) \}$
uniquely determined by the following properties:

(1) For $\pi \in NC(n)$, $k_{\pi}^{( \B )}$ is a multilinear 
functional from $\M^{n}$ to $\B$.

(2) If $\pi \in NC(p)$, $\rho \in NC(q)$, $k \in \{ 0,1, \ldots , q \}$, and
if $\sigma := \ins ( \pi \mapsto \rho ; k) \in NC(p+q)$, then for every
$x_{1}, \ldots , x_{p+q} \in \M$ we have:
\begin{equation}
\left\{  \begin{array}{l}
k_{\sigma}^{( \B )} ( x_{1}, \ldots , x_{p+q} ) \ = \ 
k_{\rho}^{( \B )} ( x_{1}, \ldots , x_{k}b, x_{k+p+1}, \ldots , x_{p+q} )  \\
                       \\
\mbox{where } b := k_{\pi}^{( \B )} ( x_{k+1}, \ldots , x_{k+p} ).
\end{array} \right.
\end{equation}

(3) For every $n \geq 1$ and $x_{1}, \ldots , x_{n} \in \M$ we have:
\begin{equation}
\sum_{\pi \in NC(n)} \ k_{\pi}^{( \B )} (x_{1}, x_{2}, \ldots , x_{n} ) 
\ = \ E( x_{1}x_{2} \cdots x_{n} ).
\end{equation}

$\ $

{\bf 5.3 Remarks and Notations.} 1) In the condition (1) of Proposition 5.2,
by ``multilinear'' we understand $\C$-multilinear. The functionals 
$k_{\pi}^{( \B )}$ turn out to actually have $\B$-multilinearity properties,
namely that:
\[
k_{\pi}^{( \B )} ( bx_{1}, x_{2}, \ldots , x_{n} ) \ = \ 
b \cdot k_{\pi}^{( \B )} ( x_{1}, x_{2}, \ldots , x_{n} ),
\]
\[
k_{\pi}^{( \B )} ( x_{1}, x_{2}, \ldots , x_{n}b ) \ = \ 
k_{\pi}^{( \B )} ( x_{1}, x_{2}, \ldots , x_{n} ) \cdot b,
\]
also that
\[
k_{\pi}^{( \B )} ( x_{1}, \ldots , x_{i}b, x_{i+1}, \ldots , x_{n} ) \ = \ 
k_{\pi}^{( \B )} ( x_{1}, \ldots , x_{i}, bx_{i+1}, \ldots , x_{n} ) 
\]
for every $\pi \in NC(n)$, $x_{1}, \ldots , x_{n} \in \M$, $b \in \B$ and 
$1 \leq i \leq n-1$. The $\C$-multilinearity stated in (1) of 
Proposition 5.2 is however more convenient when using the uniqueness part 
of the proposition.

\vspace{10pt}

2) For every $n \geq 1$, we will denote by 
$k_{n}^{( \B )}: \M^{n} \rightarrow \B$ the functional $k_{1_{n}}^{( \B )}$,
where $1_{n}$ is the partition of $\{ 1, \ldots , n \}$ into only one block.

The knowledge of the functionals $\{ k_{n}^{( \B )} \ | \ n \geq 1 \}$ 
really determines the whole family 
$\{ k_{\pi}^{( \B )} \ | \ \pi \in \cup_{n=1}^{\infty} NC(n) \}$, via the 
Equation (5.2) and the observation that the only non-crossing partitions 
that are irreducible for the operation of insertion are those of the form 
$1_{n}$. So in a certain sense the functionals $k_{\pi}^{( \B )}$ with $\pi$ 
not of the form $1_{n}$ are just some derived objects; but nevertheless, 
the $k_{\pi}^{( \B )}$'s are important for stating the essential condition 
(3) of Proposition 5.2, which can be viewed as a $\B$-valued analogue for 
Eqn.(1.11) in Section 1.

\vspace{10pt}

3) Let $( \M , E)$ and $( \widetilde{\M}, \widetilde{E} )$ be 
$\B$-probability spaces, and consider the families of elements
$x_{1}, \ldots , x_{s} \in \M$,
$\widetilde{x}_{1}, \ldots , \widetilde{x}_{s} \in \widetilde{\M}$.
We say that the families $x_{1}, \ldots , x_{s}$ and
$\widetilde{x}_{1}, \ldots , \widetilde{x}_{s}$ have 
{\em identical $\B$-cumulants} if:
\begin{equation}
k_{n}^{( \B )} (x_{r_{1}}b_{1}, \ldots , x_{r_{n-1}}b_{n-1}, x_{r_{n}} )
\ = \ 
k_{n}^{( \B )} ( \widetilde{x}_{r_{1}}b_{1}, \ldots , 
\widetilde{x}_{r_{n-1}}b_{n-1}, \widetilde{x}_{r_{n}} ),
\end{equation}
for every $n \geq 1$, $1 \leq r_{1}, \ldots , r_{n} \leq s$, and 
$b_{1}, \ldots , b_{n-1} \in \B$.

If $x_{1}, \ldots , x_{s}$ and
$\widetilde{x}_{1}, \ldots , \widetilde{x}_{s}$ have identical 
$\B$-cumulants, then the Equations (5.4) actually hold with 
``$k_{\pi}^{( \B )}$'' instead of $k_{n}^{( \B )}$; this is immediate from
(5.2), by an induction argument. Another induction argument and the use of
Equation (5.3) show that $x_{1}, \ldots , x_{s}$ and
$\widetilde{x}_{1}, \ldots , \widetilde{x}_{s}$ have identical 
$\B$-cumulants if and only if the two families are identically 
$\B$-distributed in the sense of Equation (5.1). Hence, similarly to the 
scalar case reviewed in Section 1, the $\B$-cumulants offer an alternative 
to working with $\B$-moments.

$\ $

It is useful to record the following generalization (in Proposition 5.5) of
the uniqueness part of Proposition 5.2. In all the considerations of this 
paper, by ``$\B$-bimodule'' we will understand a left-and-right $\B$-module,
where the left and the right action of $\B$ commute with each other.

$\ $

{\bf 5.4 Definition.} Let $\X$ be a $\B$-bimodule, and suppose that for every
$n \geq 1$ and $\pi \in NC(n)$ we have a $\C$-multilinear functional 
$f_{\pi} : \X^{n} \rightarrow \B$. We say that the family of functionals
$\{ f_{\pi} \ | \ \pi \in \cup_{n=1}^{\infty} NC(n) \}$ has the 
{\em insertion property} if the following holds: if
$\sigma = \ins ( \pi \mapsto \rho ; k )$ with $\pi \in NC(p)$,
$\rho \in NC(q)$, $k \in \{ 0,1, \ldots , q \}$, and if 
$x_{1}, \ldots , x_{p+q} \in \X$, then:
\begin{equation}
\left\{  \begin{array}{l}
f_{\sigma} ( x_{1}, \ldots , x_{p+q} ) \ = \ 
f_{\rho} ( x_{1}, \ldots , x_{k} \cdot b, x_{k+p+1}, \ldots , x_{p+q} ) \\
                       \\
\mbox{where } b := f_{\pi} ( x_{k+1}, \ldots , x_{k+p} ) \in \B .
\end{array} \right.
\end{equation}

$\ $

{\bf 5.5 Proposition.} Let $\X$ be a $\B$-bimodule, and suppose that for every
$n \geq 1$ and $\pi \in NC(n)$ we have two $\C$-multilinear functionals
$f_{\pi},g_{\pi} : \X^{n} \rightarrow \B$. If both the families
$\{ f_{\pi} \ | \ \pi \in \cup_{n=1}^{\infty} NC(n) \}$
and $\{ g_{\pi} \ | \ \pi \in \cup_{n=1}^{\infty} NC(n) \}$ have the 
insertion property, and if:
\begin{equation}
\sum_{\pi \in NC(n)} f_{\pi} (x_{1}, \ldots , x_{n}) \ = \ 
\sum_{\pi \in NC(n)} g_{\pi} (x_{1}, \ldots , x_{n}),
\end{equation}
for every $n \geq 1$ and $x_{1}, \ldots , x_{n} \in \X$, then we must have
that $f_{\pi} = g_{\pi}$ for all $\pi \in \cup_{n=1}^{\infty} NC(n)$.

$\ $

The proof of Proposition 5.5 is done by induction on $n$ (where 
$\pi \in NC(n)$), and is an immediate adaptation of arguments in 
\cite{S2}, Section 3.2.

The main use of $\B$-cumulants is as tool for studying freeness with 
amalgamation over $\B$. Recall that this is defined as follows (cf. e.g.
\cite{VDN}, Section 3.8).

$\ $

{\bf 5.6 Definition.} Let $( \M , E)$ be a $\B$-probability space and let 
$\M_{1}, \ldots , \M_{s}$ be subalgebras of $\M$ such that 
$\M_{1}, \ldots , \M_{s} \supset \B$. We say that $\M_{1}, \ldots , \M_{s}$ 
are {\em free with amalgamation over $\B$} if for every $n \geq 1$ and
every $r_{1}, \ldots , r_{n} \in \{ 1, \ldots , s \}$ such that 
$r_{1} \neq r_{2}, r_{2} \neq r_{3}, \ldots , r_{n-1} \neq r_{n}$ we have:
\begin{equation}
\left\{  \begin{array}{l}
x_{1} \in \M_{r_{1}}, x_{2} \in \M_{r_{2}}, \ldots , x_{n} \in \M_{r_{n}} \\
                          \\
E(x_{1}) = E( x_{2} ) = \cdots =E( x_{n} ) = 0    
\end{array}  \right\}
\ \Rightarrow \ E( x_{1}x_{2} \cdots x_{n} ) = 0.
\end{equation}

$\ $

{\bf 5.7 Remark.} The important characterization of freeness described in 
Remark 1.4 can be generalized to the $\B$-valued framework. More precisely:
if $( \M , E )$ and 
$\B \subset \M_{1} , \ldots , \M_{s} \subset \M$ are as above, then the 
freeness of $\M_{1}, \ldots , \M_{s}$ with amalgamation over $\B$ is 
equivalent to the following condition:
\[
\left\{  \begin{array}{c}
k_{n}^{( \B )} (x_{1}, \ldots , x_{n} ) = 0  \\
\mbox{whenever $x_{1} \in \M_{r_{1}}, \ldots , x_{n} \in \M_{r_{n}}$} \\
\mbox{are such that $\exists \ 1 \leq k<l \leq n$ with $r_{k} \neq r_{l}$.}
\end{array}  \right.
\]
See \cite{S2}, Section 3.3.

$\ $

{\bf 5.8 Notations.} For the remaining of this section we will suppose that 
besides the algebra $\B$ (which was fixed from the beginning of the section) 
we have also fixed:

-- a unital subalgebra $\D \subset \B$;

-- a linear map $\tau : \B \rightarrow \D$ with the properties that 
$\tau (d) = d$ for every $d \in \D$, and that 
$\tau ( d_{1}bd_{2} ) = d_{1} \tau (b) d_{2}$ for every
$d_{1},d_{2} \in \D$, $b \in \B$. 

We will assume moreover that $\tau$ is faithful (or non-degenerate) in the 
sense that if $b \in \B$ has the property that $\tau (bb') = 0$ for all 
$b' \in \B$, then $b=0$.

$\ $

In the Notations 5.8, observe that any $\B$-probability space $( \M , E)$ 
induces a $\D$-probability space $( \M , E_{\D})$, where we set
$E_{\D} := \tau \circ E$. 

$\ $

{\bf 5.9 Proposition.} Let $( \M , E)$ and $( \widetilde{\M}, \widetilde{E} )$
be $\B$-probability spaces, and consider the corresponding $\D$-probability
spaces $( \M , E_{\D} )$ and 
$( \widetilde{\M}, \widetilde{E}_{\D} )$. Suppose that ${\cal C} \subset \M$
and $\widetilde{\cal C} \subset \widetilde{\M}$ are subalgebras which contain 
$\D$, and suppose that each of ${\cal C}$ and $\widetilde{\cal C}$ is free
from $\B$ with amalgamation over $\D$ (in its corresponding space). Let 
$x_{1}, \ldots , x_{s}$ be in ${\cal C}$, and let 
$\widetilde{x}_{1}, \ldots , \widetilde{x}_{s}$ be in 
$\widetilde{\cal C}$. If the families $x_{1}, \ldots , x_{s}$ are 
identically $\D$-distributed, then the two families are also 
identically $\B$-distributed.

$\ $

{\bf Proof.} We have to show that:
\[
E_{\B} ( b_{0}x_{r_{1}}b_{1} \cdots x_{r_{n}}b_{n} ) \ = \
\widetilde{E}_{\B} ( b_{0} \widetilde{x}_{r_{1}}b_{1} \cdots 
\widetilde{x}_{r_{n}}b_{n} ) ,
\]
for every $n \geq 1$, $1 \leq r_{1}, \ldots , r_{n} \leq s$, and
$b_{0},b_{1}, \ldots , b_{n} \in \B$.
In view of the faithfulness of $\tau : \B \rightarrow \D$, this will follow 
if we can show that:
\begin{equation}
\tau ( \ E_{\B} ( b_{0}x_{r_{1}}b_{1} \cdots x_{r_{n}}b_{n} )b' \ ) \ = \
\tau ( \widetilde{E}_{\B} ( b_{0} \widetilde{x}_{r_{1}}b_{1} \cdots 
\widetilde{x}_{r_{n}}b_{n} )b' \ )  ,
\end{equation}
(for every $n, r_{1}, \ldots , r_{n}, b_{0},b_{1}, \ldots , b_{n}$ as
before, and for every $b' \in \B$). By absorbing $b'$ into $E_{\B}$ and 
into $\widetilde{E}_{\B}$, and by taking into account that 
$\tau \circ E_{\B} = E_{\D}$,
$\tau \circ \widetilde{E}_{\B} = \widetilde{E}_{\D}$, we reduce (5.8) to:
\begin{equation}
E_{\D} ( b_{0}x_{r_{1}}b_{1} \cdots x_{r_{n}}b_{n} ) \ = \
\widetilde{E}_{\D} ( b_{0} \widetilde{x}_{r_{1}}b_{1} \cdots 
\widetilde{x}_{r_{n}}b_{n} ),
\end{equation}
for every $n \geq 1$, $1 \leq r_{1}, \ldots , r_{n} \leq s$, and
$b_{0},b_{1}, \ldots , b_{n} \in \B$. Finally, (5.9) follows from the 
definition of freeness with amalgamation plus an induction argument, by 
using the hypotheses that $x_{1}, \ldots , x_{s}$ and
$\widetilde{x}_{1}, \ldots , \widetilde{x}_{s}$ have identical 
$\D$-distributions, and that ${\cal C}$, $\widetilde{\cal C}$ are free from
$\B$ with amalgamation over $\D$. {\bf QED}

$\ $

$\ $

\setcounter{section}{6}
{\large\bf 6. Cumulants with respect to the algebra of d-by-d
scalar matrices}

$\ $

{\bf 6.1 Notations.} In this section we fix a positive integer $d$, and 
\setcounter{equation}{0}
we consider the algebra $\B \ := \ M_{d} ( \C )$. If $\ncps$ is any 
non-commutative probability space, then the algebra $M_{d} (  \A )$ gets
a structure of $\B$-probability space, where we view $\B$ as a 
subalgebra of $M_{d} ( \A )$ via the natural identification:
\begin{equation}
[ \lambda_{i,j} ]_{i,j=1}^{d} \ = \ [ \lambda_{i,j} I ]_{i,j=1}^{d} 
\end{equation}
(with $I$ = the unit of $\A$). The expectation
$E_{\B}: M_{d} ( \A ) \rightarrow \B$ is defined by the formula:
\begin{equation}
E_{\B} ( \ [a_{i,j}]_{i,j=1}^{d} \ ) \ := \ 
[ \ \varphi (a_{i,j}) \ ]_{i,j=1}^{d}.
\end{equation}
Thus we are in the situation when we can consider $\B$-valued cumulants 
for families of matrices in $M_{d} ( \A )$.

The goal of the section is to give an explicit formula for 
the $\B$-valued cumulant of a family of matrices, in terms of the
scalar cumulants of the entries of these matrices.

$\ $

{\bf 6.2 Theorem.} In the framework considered above let 
$A_{1}, \ldots , A_{n}$ be matrices in $M_{d} ( \A )$, where 
$A_{m} = [ a_{i,j}^{(m)} ]_{i,j=1}^{d}$ for $1 \leq m \leq n$. Then for 
every $1 \leq i, j \leq d$, the $(i,j)$-entry $\lambda_{i,j}$ of the 
$\B$-valued cumulant $k_{n}^{( \B )} ( A_{1}, \ldots , A_{n} )$ is given 
by the formula:
\begin{equation}
\lambda_{i,j} \ = \ 
\sum_{i_{1}, \ldots , i_{n-1}=1}^{d}
k_{n} ( a_{i,i_{1}}^{(1)}, a_{i_{1},i_{2}}^{(2)}, \ldots ,
a_{i_{n-2},i_{n-1}}^{(n-1)}, a_{i_{n-1},j}^{(n)} ).
\end{equation}

$\ $

{\bf Proof.} For every $n \geq 1$ and $\pi \in NC(n)$, we define a
multilinear functional $f_{\pi}: ( M_{d}( \A ) )^{n} \rightarrow \B$,
by the following formula:
\begin{equation}
(i,j)-\mbox{entry of }
f_{\pi} (A_{1}, \ldots , A_{n}) \ := \ 
\end{equation}
\[
\sum_{i_{1}, \ldots , i_{n-1}=1}^{d} \ k_{\pi} ( a_{i,i_{1}}^{(1)},
a_{i_{1},i_{2}}^{(2)}, \ldots , 
a_{i_{n-2},i_{n-1}}^{(n-1)}, a_{i_{n-1},j}^{(n)} ) ,
\]
for every $A_{1}, \ldots , A_{n} \in M_{d} ( \A )$ and every 
$1 \leq i,j \leq d$ (and where $a_{k,l}^{(m)}$ stands for the $(k,l)$-entry
of the matrix $A_{m}$). We will verify that the family of functionals 
$\{ f_{\pi} \ | \ \pi \in \cup_{n=1}^{\infty} NC(n) \}$ satisfies the 
conditions (2) and (3) from Proposition 5.2, which determine uniquely the 
$\B$-valued cumulant functionals. Once this is done, the equality 
$f_{\pi} = k_{\pi}^{( \B )}$ applied to the partition $\pi = 1_{n}$ (of 
$\{ 1, \ldots , n \}$ into only one block) will give the statement of 
the theorem.

We start with the verification of condition (3) (about summation). Given 
$n \geq 1$ and $A_{1}, \ldots , A_{n} \in M_{d} ( \A )$, we look at:
\begin{equation}
\sum_{\pi \in NC(n)} \ f_{\pi} ( A_{1}, \ldots , A_{n} ).
\end{equation}
For every $1 \leq i,j \leq d$, the $(i,j)$-entry of the matrix appearing
in (6.5) is equal to:
\[
\sum_{\pi \in NC(n)} \ 
\sum_{i_{1}, \ldots , i_{n-1}=1}^{d} \ k_{\pi} ( a_{i,i_{1}}^{(1)},
a_{i_{1},i_{2}}^{(2)}, \ldots , 
a_{i_{n-2},i_{n-1}}^{(n-1)}, a_{i_{n-1},j}^{(n)} ) 
\]
\[
= \ \sum_{i_{1}, \ldots , i_{n-1}=1}^{d} \ 
\Bigl( \ \sum_{\pi \in NC(n)} k_{\pi} ( a_{i,i_{1}}^{(1)},
a_{i_{1},i_{2}}^{(2)}, \ldots , 
a_{i_{n-2},i_{n-1}}^{(n-1)}, a_{i_{n-1},j}^{(n)} ) \ \Bigr)
\]
\[
= \ \sum_{i_{1}, \ldots , i_{n-1}=1}^{d} \ 
\varphi ( a_{i,i_{1}}^{(1)} a_{i_{1},i_{2}}^{(2)} \cdots 
a_{i_{n-2},i_{n-1}}^{(n-1)} a_{i_{n-1},j}^{(n)} )
\]
(by the relation between scalar cumulants and moments). It is clear that 
the latter quantity is equal to $\varphi$ of the $(i,j)$-entry of 
$A_{1}A_{2} \cdots A_{n}$. Hence the matrix in (6.5) is equal to 
$E_{\B} ( A_{1}A_{2} \cdots A_{n} )$ (as desired).

We now move to the verification of condition (2) (about insertion). Suppose
that $\sigma = \ins ( \pi \mapsto \rho ; k )$, where $\pi \in NC(p)$,
$\rho \in NC(q)$, $0 \leq k \leq q$, and where the Notations 5.1 are used.
Given matrices $A_{1}, \ldots , A_{p+q} \in M_{d} ( \A )$, we want to
verify that:
\begin{equation}
f_{\sigma} (A_{1}, \ldots , A_{p+q} ) \ = \ 
f_{\rho} (A_{1}, \ldots , A_{k-1}, A_{k}B, A_{k+p+1}, \ldots , A_{p+q}), \\
\end{equation}
where 
\begin{equation}
B \ := \ f_{\pi} (A_{k+1}, \ldots , A_{k+p}) .
\end{equation}
We fix $i$ and $j$ in $\{ 1, \ldots , d \}$, and we work on the 
$(i,j)$-entry of the left-hand side of (6.6). By the definition of 
$f_{\sigma}$ this equals:
\[
\sum_{i_{1}, \ldots , i_{p+q-1}=1}^{d} \ 
k_{\sigma} ( a_{i,i_{1}}^{(1)}, a_{i_{1},i_{2}}^{(2)}, \ldots , 
a_{i_{p+q-2},i_{p+q-1}}^{(p+q-1)}, a_{i_{p+q-1},j}^{(p+q)} ),
\]
so by using the insertion property for scalar cumulants, we can re-write
it as:
\begin{equation}
\sum_{i_{1}, \ldots , i_{p+q-1}=1}^{d} \ 
k_{\pi} ( a_{i_{k},i_{k+1}}^{(k+1)}, \ldots , 
a_{i_{k+p-1},i_{k+p}}^{(k+p)} )  \cdot
\end{equation}
\[
\cdot k_{\rho} ( a_{i,i_{1}}^{(1)}, \ldots , a_{i_{k-1},i_{k}}^{(k)}, 
a_{i_{k+p},i_{k+p+1}}^{(k+p+1)}, \ldots ,
a_{i_{p+q-2},i_{p+q-1}}^{(p+q-1)}, a_{i_{p+q-1},j}^{(p+q)} ).
\]

Now, let us denote:
\begin{equation}
B =: [ \beta_{i,j} ]_{i,j=1}^{d} \in \B , \ \ 
BA_{k} =: [ x_{i,j} ]_{i,j=1}^{d} \in M_{d} ( \A ), 
\end{equation}
where $B$ is the matrix defined by (6.7). If in the summation 
of (6.8) we first sum over the indices 
$i_{k+1}, \ldots , i_{k+p-1}$, we arrive to:
\[
\sum_{i_{1}, \ldots , i_{k}, i_{k+p}, \ldots , i_{p+q-1}=1}^{d} \ 
\Bigl( \sum_{i_{k+1}, \ldots , i_{k+p-1}=1}^{d} \ 
k_{\pi} ( a_{i_{k},i_{k+1}}^{(k+1)}, \ldots , 
a_{i_{k+p-1},i_{k+p}}^{(k+p)} ) \Bigr) \cdot
\]
\[
\cdot k_{\rho} ( a_{i,i_{1}}^{(1)}, \ldots , a_{i_{k-1},i_{k}}^{(k)}, 
a_{i_{k+p},i_{k+p+1}}^{(k+p+1)}, \ldots ,
a_{i_{p+q-2},i_{p+q-1}}^{(p+q-1)}, a_{i_{p+q-1},j}^{(p+q)} )
\]
\[
= \ \sum_{i_{1}, \ldots , i_{k}, i_{k+p}, \ldots , i_{p+q-1}=1}^{d} \ 
\beta_{i_{k},i_{k+p}} \cdot
k_{\rho} ( a_{i,i_{1}}^{(1)}, \ldots , a_{i_{k-1},i_{k}}^{(k)}, 
a_{i_{k+p},i_{k+p+1}}^{(k+p+1)}, \ldots ,
a_{i_{p+q-2},i_{p+q-1}}^{(p+q-1)}, a_{i_{p+q-1},j}^{(p+q)} )
\]
(by taking into account the definition of $f_{\pi}$, the Equation (6.7),
and the notation in (6.9))
\[
= \ \sum_{i_{1}, \ldots , i_{k}, i_{k+p}, \ldots , i_{p+q-1}=1}^{d} \ 
k_{\rho} ( a_{i,i_{1}}^{(1)}, \ldots , a_{i_{k-1},i_{k}}^{(k)}
\beta_{i_{k},i_{k+p}},
a_{i_{k+p},i_{k+p+1}}^{(k+p+1)}, \ldots ,
a_{i_{p+q-2},i_{p+q-1}}^{(p+q-1)}, a_{i_{p+q-1},j}^{(p+q)} )
\]
\[
= \ \sum_{i_{1}, \ldots , i_{k-1}, i_{k+p}, \ldots , i_{p+q-1}=1}^{d} \ 
k_{\rho} ( a_{i,i_{1}}^{(1)}, \ldots , a_{i_{k-2},i_{k-1}}^{(k-1)},
x_{i_{k-1},i_{k+p}}, a_{i_{k+p},i_{k+p+1}}^{(k+p+1)}, \ldots ,
a_{i_{p+q-1},j}^{(p+q)} )
\]
(by summing over $i_{k}$ and by taking into account the definition of 
the $x_{i,j}$'s in (6.9)). The last expression is exactly the 
$(i,j)$-entry of the matrix on the right-hand side of (6.6), and 
this concludes the proof. {\bf QED}

$\ $

{\bf 6.3 Remark.} We actually arrived to prove a stronger formula than
originally announced in Theorem 6.2, namely that
\begin{equation}
(i,j)-\mbox{entry of }
k_{\pi}^{( \B )} (A_{1}, \ldots , A_{n}) \ = \ 
\end{equation}
\[
\sum_{i_{1}, \ldots , i_{n-1}=1}^{d} \ k_{\pi} ( a_{i,i_{1}}^{(1)},
a_{i_{1},i_{2}}^{(2)}, \ldots , 
a_{i_{n-2},i_{n-1}}^{(n-1)}, a_{i_{n-1},j}^{(n)} ) ,
\]
for every non-crossing partition $\pi \in NC(n)$ and every 
$A_{1}, \ldots , A_{n} \in M_{d} ( \A )$ (and where 
$a_{k,l}^{(m)}$ denotes the $(k,l)$-entry of the matrix $A_{m}$).

$\ $

{\bf 6.4 Remark.} For $A_{1}, \ldots , A_{n} \in M_{d} ( \A )$ as above,
one sometimes denotes by $A_{1} \mtens A_{2} \mtens \cdots \mtens A_{n}$ 
the matrix in $M_{d} ( \A \otimes \A \otimes  \cdots \otimes \A )$ which 
has the $(i,j)$-entry equal to:
\[
\sum_{i_{1}, \ldots , i_{n-1}=1}^{d} \ 
a_{i,i_{1}}^{(1)} \otimes a_{i_{1},i_{2}}^{(2)} \otimes \cdots
\otimes a_{i_{n-1},j}^{(n)} , \ \ 1 \leq i,j \leq d.
\]
The operation $\mtens$ is for instance used in considerations on tensor 
products of operator spaces (see e.g. Section 8.1 of \cite{ER}, or 
Section 3 of \cite{P}). 

The statement of Theorem 6.2 can be given a nice form if we use 
$\mtens$, as follows: instead of viewing the scalar-valued cumulant $k_{n}$ 
as a multilinear map from $\A^{n}$ to $\C$, let us view it as a 
linear map from the $n$-fold tensor product 
$\A \otimes \cdots \otimes \A$ into $\C$. When we go to 
$d \times d$ matrices, $k_{n}$ then induces a linear application 
$\widetilde{k_{n}}$ from $M_{d} ( \A \otimes \cdots \otimes \A )$ to 
$M_{d} ( \C )$, hence to $\B$; this is given by the formula
\[
\widetilde{k_{n}} ( \ [x_{i,j}]_{i,j=1}^{d} \ ) \ := \ 
[k_{n}(x_{i,j}) ]_{i,j=1}^{d}, \ \ 
\forall \ [ x_{i,j} ]_{i,j=1}^{d} \in M_{d} ( \A \otimes \cdots \otimes \A ).
\]
It is immediate that with these notations, the statement of Theorem 6.2 
takes the form:
\begin{equation}
k_{n}^{( \B )}( A_{1}, \ldots , A_{n} ) \ = \ 
\widetilde{k_{n}} ( A_{1} \mtens \cdots \mtens A_{n} ),
\ \ \forall \ A_{1}, \ldots , A_{n} \in M_{d} ( \A ) .
\end{equation}

$\ $

$\ $

\setcounter{section}{7}
{\large\bf 7. Cumulants with respect to the algebra of scalar diagonal
matrices}

$\ $

{\bf 7.1 Notations.} The framework for this section is similar to the one of
\setcounter{equation}{0}
Section 6, but where instead of the algebra $\B = M_{d} ( \C )$ we consider 
the algebra $\D$ of scalar diagonal $d \times d$ matrices. In other words
$\D = \mbox{span} \{ P_{1}, \ldots , P_{d} \}$, where $P_{i}$ denotes the 
matrix which has its $(i,i)$-entry equal to 1 and all the other entries
equal to 0.

If $\ncps$ is any 
non-commutative probability space, then the algebra $M_{d} (  \A )$ gets
a natural structure of $\D$-probability space, where we view $\D$ as a 
subalgebra of $M_{d} ( \A )$ via the natural identification:
\begin{equation}
\left[
\begin{array}{ccc}
\lambda_{1} &        &   0         \\
            & \ddots &             \\
0           &        & \lambda_{d} 
\end{array} \right] \ = \ 
\left[
\begin{array}{ccc}
\lambda_{1}I &        &   0         \\
             & \ddots &             \\
0            &        & \lambda_{d}I
\end{array} \right] 
\end{equation}
(with $I$ = the unit of $\A$). The expectation
$E_{\D}: M_{d} ( \A ) \rightarrow \D$ is defined by the formula:
\begin{equation}
E_{\D} ( \ [a_{i,j}]_{i,j=1}^{d} \ ) \ := \ \left[
\begin{array}{ccc}
\varphi (a_{1,1})  &        &   0         \\
                   & \ddots &             \\
      0            &        & \varphi (a_{d,d})
\end{array} \right] . 
\end{equation}
Thus we are in the situation when we can consider $\D$-valued cumulants 
for families of matrices in $M_{d} ( \A )$.

Following the same line as in the preceding section, we consider the
problem of expressing the $\D$-cumulants of a family of matrices from 
$M_{d} ( \A )$ in terms of the scalar cumulants of the entries of 
these matrices. It does not seem that there exists a nice formula holding
in general, but it is still possible to get one in the case of 
R-cyclic families. In fact we will consider a class larger than the one 
of R-cyclic families, as described in the next theorem.

$\ $

{\bf 7.2 Theorem.} In the framework considered above, let 
$A_{1}, \ldots , A_{s}$ be a family of matrices in $M_{d} ( \A )$, where 
$A_{r} = [ a_{i,j}^{(r)} ]_{i,j=1}^{d}$ for $1 \leq r \leq s$. Suppose that
for every $n \geq 1$, $1 \leq r_{1}, \ldots , r_{n} \leq s$,
$1 \leq i_{1}, \ldots , i_{n},j \leq d$ we have:
\begin{equation}
j \neq i_{n} \ \Rightarrow \ 
k_{n} ( a_{j,i_{1}}^{(r_{1})}, a_{i_{1},i_{2}}^{(r_{2})}, \ldots ,
a_{i_{n-1},i_{n}}^{(r_{n})} ) = 0.
\end{equation}
Then the $\D$-valued cumulants of the family $A_{1}, \ldots , A_{s}$ are
described by the following formula:
\begin{equation}
k_{n}^{( \D )} ( A_{r_{1}} \Lambda_{1}, \ldots , 
A_{r_{n-1}} \Lambda_{n-1}, A_{r_{n}} ) \ = 
\end{equation}
\[
\sum_{i_{1}, \ldots , i_{n}=1}^{d} \
\lambda_{i_{1}}^{(1)} \cdots \lambda_{i_{n-1}}^{(n-1)} \cdot
k_{n} ( a_{i_{n},i_{1}}^{(r_{1})}, a_{i_{1},i_{2}}^{(r_{2})}, \ldots ,
a_{i_{n-1},i_{n}}^{(r_{n})} ) P_{i_{n}},
\]
holding for $n \geq 2$, $r_{1}, \ldots , r_{n} \in \{ 1, \ldots , s \}$,
and where
\[
\Lambda_{k} \ := \ \left[
\begin{array}{ccc}
\lambda_{1}^{(k)}     &        &   0         \\
                      & \ddots &             \\
         0            &        & \lambda_{d}^{(k)}
\end{array} \right] \in \D , \ \ 1 \leq k \leq n-1.
\]

$\ $

{\bf Proof.} Let $\X$ be the free $\D$-bimodule with $s$ generators
$X_{1}, \ldots , X_{s}$. As a vector space over $\C$, $\X$ has dimension
$d^{2}s$, and has a natural basis given by the elements
$P_{i}X_{r}P_{j}$, with $1 \leq i,j \leq d$ and $1 \leq r \leq s$.

For every $n \geq 1$ and $\pi \in NC(n)$ we consider the 
$\C$-multilinear functionals $f_{\pi}$ and $g_{\pi}$ from $\X^{n}$ to 
$\D$, determined as follows (by their action on the natural basis
of $\X^{n}$):
\begin{equation}
f_{\pi} (P_{i_{1}}X_{r_{1}}P_{j_{1}} , \ldots ,
P_{i_{n}}X_{r_{n}}P_{j_{n}} ) \ = \ 
k_{\pi}^{( \D )} (P_{i_{1}}A_{r_{1}}P_{j_{1}} , \ldots ,
P_{i_{n}}A_{r_{n}}P_{j_{n}} ) 
\end{equation}
and
\begin{equation}
g_{\pi} (P_{i_{1}}X_{r_{1}}P_{j_{1}} , \ldots ,
P_{i_{n}}X_{r_{n}}P_{j_{n}} ) \ = \ 
\end{equation}
\[
\delta_{i_{1},j_{n}} \delta_{i_{2},j_{1}} \cdots 
\delta_{i_{n},j_{n-1}} \cdot k_{\pi} ( a_{j_{n},j_{1}}^{(r_{1})},
a_{j_{1},j_{2}}^{(r_{2})}, \ldots ,
a_{j_{n-1},j_{n}}^{(r_{n})} ) P_{j_{n}},
\]
for $n \geq 1$ and $1 \leq r_{1}, \ldots , r_{n} \leq s$,
$1 \leq i_{1}, \ldots , i_{n},j_{1}, \ldots , j_{n} \leq d$.
An immediate linearity argument shows that:
\begin{equation}
f_{\pi} ( \Gamma_{1} X_{r_{1}} \Lambda_{1} , \ldots ,
\Gamma_{n} X_{r_{n}} \Lambda_{n} ) \ = \ 
k_{\pi}^{( \D )} ( \Gamma_{1} A_{r_{1}} \Lambda_{1} , \ldots ,
\Gamma_{n} A_{r_{n}} \Lambda_{n} ), 
\end{equation}
and
\begin{equation}
g_{\pi} ( \Gamma_{1} X_{r_{1}} \Lambda_{1} , \ldots ,
\Gamma_{n} X_{r_{n}} \Lambda_{n} ) \ = \ 
\end{equation}
\[
\sum_{j_{1}, \ldots , j_{n}=1}^{d} \ \gamma_{j_{n}}^{(1)} \cdot
\lambda_{j_{1}}^{(1)} \gamma_{j_{1}}^{(2)} \cdots
\lambda_{j_{n-1}}^{(n-1)} \gamma_{j_{n-1}}^{(n)} \cdot
\lambda_{j_{n}}^{(n)} \cdot
k_{\pi} ( a_{j_{n},j_{1}}^{(r_{1})}, a_{j_{1},j_{2}}^{(r_{2})},
\ldots , a_{j_{n-1},j_{n}}^{(r_{n})} ) P_{j_{n}},
\]
for every $n \geq 1$, $1 \leq r_{1}, \ldots , r_{n} \leq s$ and 
$\Gamma_{1}, \Lambda_{1}, \ldots , \Gamma_{n}, \Lambda_{n} \in \D$,
where:
\[
\Gamma_{k} \ := \ \left[
\begin{array}{ccc}
\gamma_{1}^{(k)}    &        &   0         \\
                    & \ddots &             \\
         0          &        & \gamma_{d}^{(k)}
\end{array} \right] , \ \ 
\Lambda_{k} \ := \ \left[
\begin{array}{ccc}
\lambda_{1}^{(k)}  &        &   0         \\
                   & \ddots &             \\
         0         &        & \lambda_{d}^{(k)}
\end{array} \right] , \ \ 1 \leq k \leq n.
\]

Now, both the families of functionals
$\{ f_{\pi} \ | \ \pi \in \cup_{n=1}^{\infty} NC(n) \}$ and
$\{ g_{\pi} \ | \ \pi \in \cup_{n=1}^{\infty} NC(n) \}$ satisfy the insertion
property considered in Definition 5.4. For the $f_{\pi}$'s this is an 
immediate consequence of the corresponding property for the $\D$-valued
cumulant functionals
$\{ k_{\pi}^{( \D )} \ | \ \pi \in \cup_{n=1}^{\infty} NC(n) \}$.
For the $g_{\pi}$'s the insertion property follows from a calculation very
similar in nature to the one shown in the proof of Theorem 6.2, and which,
due to its routine character, will be left to the reader. (The reader who 
will have the patience to go through this calculation will notice that 
it effectively makes use of the implication (7.3) stated in the 
hypothesis of the current theorem.)

We next show that $f_{\pi} = g_{\pi}$ for every 
$\pi \in \cup_{n=1}^{\infty} NC(n)$. Proposition 5.5 combined with a 
linearity argument shows that all we need to check is the equality:
\[
\sum_{\pi \in NC(n)} g_{\pi} (P_{i_{1}}X_{r_{1}}P_{j_{1}} , \ldots ,
P_{i_{n}}X_{r_{n}}P_{j_{n}} ) \ = \ \sum_{\pi \in NC(n)} f_{\pi}
(P_{i_{1}}X_{r_{1}}P_{j_{1}} , \ldots , P_{i_{n}}X_{r_{n}}P_{j_{n}} )
\]
(for some fixed $n \geq 1$, $1 \leq r_{1}, \ldots , r_{n} \leq s$,
$1 \leq i_{1}, j_{1}, \ldots , i_{n}, j_{n} \leq d$). And indeed, 
we compute:
\[
\sum_{\pi \in NC(n)} g_{\pi} (P_{i_{1}}X_{r_{1}}P_{j_{1}} , \ldots ,
P_{i_{n}}X_{r_{n}}P_{j_{n}} ) \ = \ 
\]
\[
\delta_{i_{1},j_{n}} \delta_{i_{2},j_{1}} \cdots 
\delta_{i_{n},j_{n-1}} \cdot \Bigl( \sum_{\pi \in NC(n)}
k_{\pi} ( a_{j_{n},j_{1}}^{(r_{1})}, a_{j_{1},j_{2}}^{(r_{2})}, 
\ldots , a_{j_{n-1},j_{n}}^{(r_{n})} ) \  \Bigr) P_{j_{n}}
\]
\[
= \ \delta_{i_{1},j_{n}} \delta_{i_{2},j_{1}} \cdots 
\delta_{i_{n},j_{n-1}} \cdot \varphi (
a_{j_{n},j_{1}}^{(r_{1})} a_{j_{1},j_{2}}^{(r_{2})} 
\cdots , a_{j_{n-1},j_{n}}^{(r_{n})} ) P_{j_{n}}
\]
\[
= \ E_{\D} (P_{i_{1}}A_{r_{1}}P_{j_{1}} 
P_{i_{2}}A_{r_{2}}P_{j_{2}} \cdots  P_{i_{n}}A_{r_{n}}P_{j_{n}} ) 
\]
\[
= \ \sum_{\pi \in NC(n)} k_{\pi}^{( \D )}
(P_{i_{1}}A_{r_{1}}P_{j_{1}} , \ldots  , P_{i_{n}}A_{r_{n}}P_{j_{n}} ) 
\]
\[
= \ \sum_{\pi \in NC(n)} f_{\pi}
(P_{i_{1}}X_{r_{1}}P_{j_{1}} , \ldots  , P_{i_{n}}X_{r_{n}}P_{j_{n}} ).
\]

But if $f_{\pi} = g_{\pi}$, then one can equate the right-hand sides of the 
Equations (7.7) and (7.8). By doing this for $\pi = 1_{n}$ (the partition
of $\{ 1, \ldots , n \}$ into only one block), and by appropriately
choosing $\Gamma_{1}, \Lambda_{1}, \ldots , \Gamma_{n},
\Lambda_{n} \in \D$, one obtains the Equation (7.4) from the conclusion 
of the theorem. {\bf QED}

$\ $

{\bf 7.3 Remark.} In the framework of the Notations 7.1, let 
$A_{1}, \ldots , A_{s}$ be an R-cyclic family of matrices in $M_{d} ( \A )$,
where $A_{r} = [ a_{i,j}^{(r)} ]_{i,j=1}^{d}$ for $1 \leq r \leq s$. Then
Theorem 7.2 gives us an interpretation for the cyclic cumulants of the 
entries of $A_{1}, \ldots , A_{s}$ (i.e., for the coefficients of the 
determining series of the family $A_{1}, \ldots , A_{s}$). More precisely,
for every $n \geq 1$, $1 \leq r_{1}, \ldots , r_{n} \leq s$ and 
$1 \leq i_{1}, \ldots , i_{n} \leq d$, we have that:
\begin{equation}
k_{n} ( \ a_{i_{n},i_{1}}^{(r_{1})}, a_{i_{1},i_{2}}^{(r_{2})}, \ldots ,
a_{i_{n-2},i_{n-1}}^{(r_{n-1})}, a_{i_{n-1},i_{n}}^{(r_{n})}  \ ) \ =  
\end{equation}
\[
(i_{n},i_{n})-\mbox{entry of }
k_{n}^{( \D )} ( A_{r_{1}}P_{i_{1}}, \ldots , 
A_{r_{n-1}} P_{i_{n-1}} ,  A_{r_{n}} ) .
\]

$\ $

{\bf 7.4 Remark.} In analogy to Remark 6.4, one can also reformulate 
the result of Theorem 7.2 by using the $\mtens$-product. Let us denote 
by $\widetilde{k_n}^\D$ the counterpart of $\widetilde{k_{n}}$ (from
Remark 6.4) which is suitable for working with $\D$. That is,
$\widetilde{k_{n}}^\D$ is the linear application from 
$M_{d} ( \A \otimes \cdots \otimes \A )$ to
$\D$ given by the formula:
\[
\widetilde{k_{n}}^\D ( \ [x_{i,j}]_{i,j=1}^{d} \ ) \ := \
[ \delta_{i,j} k_{n}(x_{i,j}) ]_{i,j=1}^{d} , \ \ \forall \
[ x_{i,j} ]_{i,j=1}^{d} \in M_{d} ( \A \otimes \cdots \otimes \A ).
\]
It is immediate that with these notations, the statement of Theorem 
7.2 takes the following form (where we have set the matrices 
$\Lambda_{1}, \ldots , \Lambda_{n-1}$ from Equation (7.4) to be equal 
to the unit of $\D$):
\begin{equation}
k_{n}^{( \D )}( A_{r_1}, \ldots , A_{r_n} ) \ = \
\widetilde{k_{n}}^\D ( A_{r_1} \mtens \cdots \mtens A_{r_n} ),
\ \ \forall \ 1\leq r_1,\dots,r_n\leq s .
\end{equation}

It is a natural question if the same kind of formula is true when 
we consider other algebras of scalar $d \times d$ matrices (instead of 
$\B$ and $\D$, as we have in the Equations (6.11) and (7.10)). Let us 
consider the case of the smallest possible such algebra, namely $\C$
(corresponding to scalar multiples of the identity $d \times d$ matrix).
Again we have a linear application
$\widetilde{k_{n}}^\C : M_{d} ( \A \otimes \cdots \otimes \A ) 
\rightarrow \C$, given by the formula:
\[
\widetilde{k_{n}}^\C ( \ [x_{i,j}]_{i,j=1}^{d} \ ) \ := \
\frac 1d \sum_{i=1}^d k_{n}(x_{i,i}) , \ \ \forall \
[ x_{i,j} ]_{i,j=1}^{d} \in M_{d} ( \A \otimes \cdots \otimes \A ).
\]
The question becomes: under what conditions on the matrices 
$A_{1}, \ldots , A_{s} \in M_{d} ( \A )$ can we infer that:
\begin{equation}
k_{n} ( A_{r_1}, \ldots , A_{r_n} ) \ = \
\widetilde{k_{n}}^\C ( A_{r_1} \mtens \cdots \mtens A_{r_n} ),
\end{equation}
for every $1\leq r_1,\dots,r_n\leq s$?
It turns out (see \cite{NSS4}) that (7.11) can be guaranteed if we know that
$k_{n}^{( \D )}( A_{r_1}, \ldots , A_{r_n} )$, which is a priori an element
in $\D$, is actually an element in $\C$. In the context of R-cyclic 
matrices, this amounts precisely to the situation discussed in Proposition 
3.1; indeed, the ``partial summation property'' stated in Equation (3.1) 
asks that $k_{n}^{( \D )}( A_{r_1}, \ldots , A_{r_n} )$ is a 
scalar multiple of the identity $d \times d$ matrix, while on the other 
hand the conclusion of Proposition 3.1 (as appearing e.g. in Equation (3.3))
is tantamount to (7.11).

$\ $

$\ $

\setcounter{section}{8}
{\large\bf 8. Characterization of R-cyclicity as freeness with amalgamation}

$\ $

In this section we combine the frameworks used in the Sections 6 and 7.
\setcounter{equation}{0}
That is, for a fixed integer $d \geq 1$ we will consider both the algebra 
$\B = M_{d} ( \C )$ and its subalgebra $\D$ consisting of diagonal matrices.
For every $1 \leq i,j \leq d$ we will denote by $V_{i,j} \in \B$ the matrix 
which has 1 on the $(i,j)$-entry and 0 on all the other entries. (Note that
the matrices denoted up to now by ``$P_{i}$'' have become ``$V_{i,i}$'', 
for $1 \leq i \leq d$.)

If $\ncps$ is any non-commutative probability space, then $M_{d} (  \A )$ 
is at the same time a $\B$-probability space and a $\D$-probability space,
where the identifications $\D \subset \B \subset M_{d} ( \A)$ and the 
expectations $E_{\B}: M_{d} ( \A ) \rightarrow \B$,
$E_{\D}: M_{d} ( \A ) \rightarrow \D$ are as described in the Sections 6
and 7. Note that the restriction of $E_{\D}$ to $\B$ is 
faithful; this implies that the discussion concluding the Section 5
(and in particular the Proposition 5.9) can be applied in this framework.

$\ $

{\bf 8.1 Lemma.} In the framework considered above, let 
$C_{1}, C_{2}, \ldots , C_{n} \in M_{d} ( \A )$ form an R-cyclic family,
where $n \geq 2$. Suppose that:

(i) for $m \in \{ 1, n \}$ we have that either $E_{\D} (C_{m}) = 0$ or 
that $C_{m}$ is the unit of $M_{d} ( \A )$; and 

(ii) for $m \in \{ 2,3, \ldots ,n-1 \}$ we have that $E_{\D} (C_{m}) = 0$.
\newline
Consider also some indices 
$i_{1}, j_{1}, \ldots , i_{n-1},j_{n-1} \in \{ 1, \ldots , d \}$ such that 
$i_{1} \neq j_{1}, \ldots , i_{n-1} \neq j_{n-1}$. Then:
\begin{equation}
E_{\D} ( \  C_{1}V_{i_{1},j_{1}} \cdots
C_{n-1} V_{i_{n-1},j_{n-1}}C_{n} \ ) \ = \ 0. 
\end{equation}

$\ $

{\bf Proof.} We will denote by $c_{i,j}^{(m)}$ the $(i,j)$-entry of $C_{m}$
($1 \leq i,j \leq d$, $1 \leq m \leq n$). The hypotheses (i) and (ii) given
above show that:
\begin{equation}
\left\{  \begin{array}{lcl}
\varphi ( c_{i,i}^{(m)} ) = 0  \mbox{ or } c_{i,i}^{(m)} = I  &
\mbox{if}  & m \in \{ 1, n \} , \  1 \leq i \leq d                   \\
\varphi ( c_{i,i}^{(m)} ) = 0   &
\mbox{if}  & m \in \{ 2,3, \ldots ,n-1 \} , \ 1 \leq i \leq d.
\end{array}  \right.
\end{equation}
In connection to this, let us also record the fact that:
\begin{equation}
\varphi ( c_{i,j}^{(m)} ) = 0, \  \ \forall \ 1 \leq m \leq n, 
\forall \ 1 \leq i,j \leq d \ \mbox{such that } i \neq j,
\end{equation}
which follows from R-cyclicity 
($\varphi ( c_{i,j}^{(m)} ) = k_{1} ( c_{i,j}^{(m)} ) = 0$ for $i \neq j$).

We will present the proof under the assumption that $n \geq 3$. The
(similar, and simpler) case $n=2$ is left as an exercise to the reader.

If we write explicitly the $(i,i)$-entry of the scalar diagonal matrix on
the left-hand side of (8.1), it becomes clear that what we have to do 
in this proof is to fix an $i \in \{ 1, \ldots , d \}$, and show that
\begin{equation}
\varphi ( \ c_{i,i_{1}}^{(1)} c_{j_{1},i_{2}}^{(2)} \cdots 
c_{j_{n-2},i_{n-1}}^{(n-1)} c_{j_{n-1},i}^{(n)} \ ) \ = \ 0.
\end{equation}
By using the relation between moments and non-crossing cumulants, the
quantity on the left-hand side of (8.4) can be written as:
\begin{equation}
\sum_{\pi \in NC(n)} \ k_{\pi} (
\ c_{i,i_{1}}^{(1)}, c_{j_{1},i_{2}}^{(2)}, \ldots , 
c_{j_{n-2},i_{n-1}}^{(n-1)}, c_{j_{n-1},i}^{(n)} \ ) .
\end{equation}
We will actually prove that every term of the sum in (8.5) is 
equal to 0.

So, besides $i \in \{ 1, \ldots , d \}$, let us also fix a partition
$\pi \in NC(n)$, and let us examine the non-crossing cumulant
$k_{\pi} ( \ c_{i,i_{1}}^{(1)}, c_{j_{1},i_{2}}^{(2)}, \ldots , 
c_{j_{n-2},i_{n-1}}^{(n-1)}, c_{j_{n-1},i}^{(n)} \ )$. Recall from
Section 1.3 that this cumulant is defined as a product having as
many factors as there are blocks in $\pi$. For the sake of brevity,
we will denote it in the rest of the proof by just 
``$k_{\pi}$''.

Denoting by $B$ the block 
of $\pi$ which contains the number 2, we distinguish four cases:

\vspace{10pt}

{\em Case 1.} $B = \{ 2 \}$.
In this case, $k_{\pi}$ 
has a factor ``$k_{1} ( c_{j_{1},i_{2}}^{(2)} )$'', which is equal to 
0 by (8.2), (8.3). So $k_{\pi}$ itself is equal to 0.

\vspace{10pt}

{\em Case 2.} $B = \{ 1, 2 \}$.
In this case, $k_{\pi}$ has a factor
``$k_{2} ( c_{i,i_{1}}^{(1)}, c_{j_{1},i_{2}}^{(2)} )$'', which is equal to
0 by R-cyclicity and the hypothesis that $i_{1} \neq j_{1}$. So again
$k_{\pi} = 0$.

\vspace{10pt}

{\em Case 3.} $B \ni 3$.
In this case, $k_{\pi}$ has a factor 
``$k_{|B|} ( \ldots, c_{j_{1},i_{2}}^{(2)}, c_{j_{2}, i_{3}}^{(3)}, 
\ldots )$'', which is equal to 0 by R-cyclicity and the hypothesis that 
$i_{2} \neq j_{2}$. So again $k_{\pi} = 0$.

\vspace{10pt}

{\em Case 4.} $B$ does not fall in any of the Cases 1-3.
In this case $B$ intersects $\{ 4, \ldots , n \}$; let us denote
$m := \min ( B \cap \{ 4, \ldots , n \} )$. The set 
$\{ 3, 4, \ldots , m-1 \}$ is a union of blocks of $\pi$; because $\pi$ 
is non-crossing, there is one of these blocks, $B_{1}$, which has to 
be an interval-block  -- say that $B_{1} = [ p, q ] \cap \Z$, with
$3 \leq p \leq q \leq m-1 \ ( \leq n-1 )$. The cumulant $k_{\pi}$ has a 
factor $k_{|B_{1}|} ( \cdots )$ corresponding to the block $B_{1}$. If
$B_{1}$ has only one element (i.e. $p=q$), then the factor 
$k_{|B_{1}|} ( \cdots )$ is equal to 0 by the same argument as in Case 1;
while if $|B_{1}| > 1$ (i.e. $p<q$), then the factor
$k_{|B_{1}|} ( \cdots )$ is equal to 0 by the same argument as 
in Cases 2, 3. Either way, $k_{\pi}$ is equal to 0. {\bf QED}

$\ $

{\bf 8.2 Theorem.} Let $A_{1}, \ldots , A_{s}$ be a family of 
matrices in $M_{d} ( \A )$, and let ${\cal C}$ denote the subalgebra
of $M_{d} ( \A )$ generated by $\{ A_{1}, \ldots , A_{s} \} \cup \D$.
The family $A_{1}, \ldots , A_{s}$ is R-cyclic if and only if 
${\cal C}$ is free from $\B$, with amalgamation over $\D$.

$\ $

{\bf Proof.} ``$\Rightarrow$''. We will verify that ${\cal C}$ is free
from $\B$, with amalgamation over $\D$, by using the definition of 
freeness with amalgamation. That is: we consider an alternating sequence
$X_{1}, X_{2}, \ldots , X_{k}$ of matrices from $\B$ and from ${\cal C}$,
such that $E_{\D} (X_{1}) = \cdots = E_{\D} ( X_{k} ) = 0$, and we want 
to show that $E_{\D} ( X_{1}X_{2} \cdots X_{k}) = 0$.

If the alternating sequence of matrices considered above does not begin
with a matrix from ${\cal C}$, let us add on the left end of the sequence 
one more matrix, equal to the identity of $M_{d} ( \A )$, and viewed as
belonging to ${\cal C}$. Let us also use this procedure at the right end
of the alternating sequence. With these adjustments we can assume that 
$k$ (the number of matrices in the sequence) is odd, $k =2n-1$, and that 
$X_{1}, X_{2n-1} \in {\cal C}$. On the other hand the hypothesis which 
we have on $X_{1}$ and $X_{2n-1}$ has to be weakened to the fact that 
they either have zero $\D$-expectation, or they are equal to the 
identity of $M_{d} ( \A )$.

We re-denote the matrices $X_{1},X_{3}, \ldots , X_{2n-1}$ by 
$C_{1}, \ldots , C_{n} ( \in {\cal C} )$. The family 
$C_{1}, \ldots , C_{n}$ is R-cyclic, by Theorem 4.3.

On the other hand, let us look at the matrices $X_{2},X_{4}, \ldots ,
X_{2n-2}$, which belong to $\B$ and have $\D$-expectation equal to 0.
It is clear that each of these matrices belongs to 
$\mbox{span} \{ V_{i,j} \ | \ 1 \leq i,j \leq d, \ i \neq j \}$.
An immediate argument with linear combinations allows us to assume 
without loss of generality that in fact we have
$X_{2} = V_{i_{1},j_{1}}, \ldots , X_{2n-2} = V_{i_{n-1},j_{n-1}}$
for some $i_{1},j_{1}, \ldots , i_{n-1},j_{n-1} \in \{ 1, \ldots , d \}$
such that $i_{1} \neq j_{1}, \ldots , i_{n-1} \neq j_{n-1}$.

With the above adjustments, the product $X_{1}X_{2} \cdots X_{k}$ now
reads: $C_{1}V_{i_{1},j_{1}} \cdots C_{n-1}$
$V_{i_{n-1},j_{n-1}}C_{n}$. The fact that this product has
zero $\D$-expectation is exactly what was proved in Lemma 8.1. 

\vspace{10pt}

``$\Leftarrow$'' In a different non-commutative probability space 
$( \N , \psi )$ we construct a family of elements 
$\{ x_{i,j}^{(r)} \ | \ 1 \leq i,j \leq d, \ 1 \leq r \leq s \}$ such that:
\begin{equation}
k_{n} ( x_{i_{1},j_{1}}^{(r_{1})}, \ldots , x_{i_{n},j_{n}}^{(r_{n})} ) 
\ = \ 0
\end{equation}
for every $n \geq 1$ and $1 \leq r_{1}, \ldots , r_{n} \leq s$,
$1 \leq i_{1},j_{1}, \ldots , i_{n},j_{n} \leq d$ for which it is not true
that $j_{1}=i_{2}, \ldots , j_{n-1}=i_{n}, j_{n}=i_{1}$; and such that
\begin{equation}
k_{n} ( x_{i_{n},i_{1}}^{(r_{1})}, x_{i_{1},i_{2}}^{(r_{2})} , \ldots ,  
x_{i_{n-1},i_{n}}^{(r_{n})} ) \ = 
\end{equation}
\[
(i_{n},i_{n})-\mbox{entry of }
k_{n}^{( \D )} ( A_{r_{1}}V_{i_{1},i_{1}}, \ldots ,
A_{r_{n-1}}V_{i_{n-1},i_{n-1}}, A_{r_{n}} ),
\]
for every $n \geq 1$ and $1 \leq r_{1}, \ldots , r_{n} \leq s$,
$1 \leq i_{1}, \ldots , i_{n} \leq d$.
Such a construction is possible because one can in general construct families
of elements with any prescribed family of scalar cumulants, via an abstract
free product construction (see e.g. \cite{VDN}, Chapter 1).

For $1 \leq r \leq s$, we consider the matrix
$X_{r} = [ x_{i,j}^{(r)} ]_{i,j=1}^{d} \in M_{d} ( \N )$. The 
Equations (8.6) and (8.7) tell us that the family $X_{1}, \ldots , X_{s}$
is R-cyclic.

Observe that for every $n \geq 1$ and every 
$1 \leq r_{1}, \ldots , r_{n} \leq s$, $1 \leq i_{1}, \ldots ,i_{n} \leq d$
we have that:
\begin{equation}
k_{n}^{( \D )} ( X_{r_{1}}V_{i_{1},i_{1}}, \ldots ,
X_{r_{n-1}}V_{i_{n-1},i_{n-1}}, X_{r_{n}} ) \ = \ 
\end{equation}
\[
k_{n}^{( \D )} ( A_{r_{1}}V_{i_{1},i_{1}}, \ldots ,
A_{r_{n-1}}V_{i_{n-1},i_{n-1}}, A_{r_{n}} ).
\]
Indeed, the scalar diagonal matrices on both sides of the Equation (8.8) have 
their $(j,j)$-entry equal to the cumulant
$k_{n} ( x_{j,i_{1}}^{(r_{1})}, x_{i_{1},i_{2}}^{(r_{2})} , \ldots ,  
x_{i_{n-1},j}^{(r_{n})} )$ (for the right-hand side this is just (8.7),
while for the left-hand side we invoke Equation (7.9) from Remark 7.3). By 
taking linear combinations with respect to $V_{i_{1},i_{1}}, \ldots , 
V_{i_{n-1},i_{n-1}}$ in (8.8) we find that the family $X_{1}, \ldots , X_{s}$
has identical $\D$-cumulants with the family $A_{1}, \ldots , A_{s}$.
Or in other words, the families $A_{1}, \ldots , A_{s}$ and
$X_{1}, \ldots , X_{s}$ have identical $\D$-distributions.

Now, our current hypothesis is that the algebra ${\cal C}$ generated by 
$\{ A_{1}, \ldots , A_{s} \} \cup \D$ is free from $\B$, with amalgamation 
over $\D$. On the other hand, the same is true about the algebra 
$\widetilde{\cal C} \subset \N$ generated by 
$\{ X_{1}, \ldots , X_{s} \} \cup \D$; this follows from the fact that the 
family $X_{1}, \ldots , X_{s}$ is R-cyclic, and the implication 
``$\Rightarrow$'' (proved above!) of the current theorem. But then we are 
in the position to apply the Proposition 5.9, which gives us that the 
families $A_{1}, \ldots , A_{s}$ and $X_{1}, \ldots , X_{s}$ actually have 
identical $\B$-distributions. The latter fact implies in turn that
the families $A_{1}, \ldots , A_{s}$ and $X_{1}, \ldots , X_{s}$ have 
identical $\B$-cumulants.

Finally, let us fix $n \geq 1$, $r_{1}, \ldots , r_{n} \in \{ 1, \ldots ,s \}$
and $i_{1},j_{1}, \ldots , i_{n},j_{n} \in \{ 1, \ldots , d \}$, and suppose 
it is not true that $j_{1}=i_{2}, \ldots , j_{n-1}=i_{n}, j_{n}=i_{1}$.
Then:
\[
k_{n} ( a_{i_{1},j_{1}}^{(r_{1})}, \ldots , a_{i_{n},j_{n}}^{(r_{n})} ) 
\]
\[
= \ (i_{1},j_{n})-\mbox{entry of }
k_{n}^{( \B )} ( A_{r_{1}}V_{j_{1},i_{2}}, \ldots ,
A_{r_{n-1}}V_{j_{n-1},i_{n}}, A_{r_{n}} ) \mbox{ (by Theorem 6.2)}
\]
\[
= \ (i_{1},j_{n})-\mbox{entry of }
k_{n}^{( \B )} ( X_{r_{1}}V_{j_{1},i_{2}}, \ldots ,
X_{r_{n-1}}V_{j_{n-1},i_{n}}, X_{r_{n}} )
\]
(since the families $A_{1}, \ldots , A_{s}$ and $X_{1}, \ldots , X_{s}$
have identical $\B$-cumulants)
\[
= \ k_{n} ( x_{i_{1},j_{1}}^{(r_{1})}, \ldots , x_{i_{n},j_{n}}^{(r_{n})} ) 
\ \mbox{ (again by Theorem 6.2)}
\]
\begin{center}
$= \ 0$ (by Equation (8.6)).   {\bf QED}
\end{center}

$\ $

$\ $

\end{document}